\documentclass[a4paper,11pt]{article}
\usepackage{float} 
\usepackage{fancyhdr}
\usepackage{latexsym,amssymb,amsmath}
\usepackage{amsmath}
\usepackage[french]{babel}
\usepackage[latin1]{inputenc}

\input xy
\xyoption{all}

\newtheorem{defi}{D\'efinition}[section]
\newtheorem{thm}[defi]{Th\'eor\`eme}
\newtheorem{prop}[defi]{Proposition}
\newtheorem{coro}[defi]{Corollaire}
\newtheorem{lemme}[defi]{Lemme}
\newtheorem{rema}[defi]{Remarque}

\newtheorem{conv}[defi]{Convention}

\newtheorem{nota}[defi]{Notation}
\newtheorem{notas}[defi]{Notations}

\newenvironment{dem}{\noindent {\it D\'emonstration $-$ }
                 \noindent}{\hfill $\Box$\vskip 5mm}

\newcommand{\N}{\mathbb{N}}
\newcommand{\Z}{\mathbb{Z}}
\newcommand{\Q}{\mathbb{Q}}
\newcommand{\R}{\mathbb{R}}
\newcommand{\C}{\mathbb{C}}

\newcommand{\PP}{\mathcal{P}}

\newcommand{\hh}{\mathcal{H}}

\newcommand{\Sym}{\mbox{Sym}}

\newcommand{\Hom}{\mbox{Hom}}
\newcommand{\Ext}{\mbox{Ext}}
\newcommand{\Id}{\mbox{Id}}
\newcommand{\isom}{\overset{\sim}{\to}}
 
\newcommand{\OO}{\mathcal{O}}  
\newcommand{\cinf}{\mathcal{C}^{\infty}} 

\newcommand{\For}{\mbox{For}}

\newcommand{\G}{\mathbb{G}}
\newcommand{\Res}{\mbox{Res}}
\newcommand{\Eis}{\mathcal{E}is}
\newcommand{\Tr}{\mbox{Tr}}
\newcommand{\Norm}{\mbox{N}}

\newcommand{\h}{\mathfrak{H}}
\newcommand{\Sp}{\mathbb{S}}
\newcommand{\X}{\mathfrak{X}}

\begin{document}

\hfill{9 Février 2008}
\begin{center}
\vspace{1cm}
{\huge  {\bf Les classes d'Eisenstein des vari\'et\'es }}
 \\ 
 \vspace{0.3cm}
 {\huge {\bf de Hilbert-Blumenthal}} \\
\vspace{1cm}
{\Large David Blotti\`ere}\\
$\;$\\
Universit\"at Paderborn,\\
Institut f\"ur Mathematik, \\
Warburger Str. 100,\\
33098 Paderborn, Germany.\\
email: \verb+blottier@math.upb.de+ \\
\vspace{1cm}

\noindent {\Large {\bf R\'esum\'e}} \\

\begin{tabular}{p{15cm}}
\\
Dans \cite{b}, on a d\'emontr\'e que les courants de Levin (cf. \cite{l}) 
permettent de d\'ecrire explicitement
les classes d'Eisenstein d'un sch\'ema ab\'elien au niveau topologique. 
On applique ici ce r\'esultat, conjectur\'e par Levin,
au cas o\`u le sch\'ema ab\'elien est une famille de vari\'et\'es ab\'eliennes de Hilbert-Blumenthal
(cf. Proposition \ref{explicit_Eis}). 
On \'etudie ensuite la d\'eg\'en\'erescence de ces classes d'Eisenstein
en une pointe de la compactification de Baily-Borel de la vari\'et\'e de Hilbert-Blumenthal. Au moyen du Th\'eor\`eme
de Burgos-Wildeshaus \cite[Theorem 2.9]{bw}, on d\'emontre un r\'esultat de rigidit\'e 
(cf. Proposition \ref{prig}) qui permet
de restreindre l'\'etude au niveau topologique. On prouve, en utilisant la description explicite des classes d'Eisenstein obtenue pr\'ec\'edemment, que ces classes d\'eg\'en\`erent en des valeurs sp\'eciales d'une fonction $L$ associ\'ee au corps de nombres totalement r\'eel sous-jacent 
(Th\'eor\`eme \ref{R1}). 
On en d\'eduit une preuve g\'eom\'etrique du Th\'eor\`eme de
Klingen-Siegel (Corollaire \ref{R2}) et un r\'esultat de non annulation pour certaines de ces classes d'Eisenstein (Corollaire \ref{R3}).\\ \\ \\
\end{tabular}

\noindent {\Large {\bf Abstract}} \\

\begin{tabular}{p{15cm}}
 \\
In \cite{b}, we have proved that Levin's currents (cf. \cite{l}) 
give an explicit description of the Eisenstein
classes of an abelian scheme at the topological level. 
We apply here this result, conjectured by Levin, in the
situation where the abelian scheme is an Hilbert-Blumenthal family of abelian varieties 
(cf. Proposition \ref{explicit_Eis}). Then we study the degeneration
of these Eisenstein classes at a cusp of the Baily-Borel compactification of the Hilbert-Blumenthal variety.
Using the Theorem of Burgos-Wildeshaus \cite[Theorem 2.9]{bw}, we prove a rigidity result which allows us to
restrict the study at the topological level (cf. Proposition \ref{prig}). 
We show, using the explicit description of the Eisenstein classes obtained previously, that these classes degenerate in special
values of an $L$-function associated to the underlying totally
real number field (Th\'eor\`eme \ref{R1}). We deduce then a geometric proof the Klingen-Siegel 
Theorem (Corollaire \ref{R2}) and
a non vanishing result for some of these Eisenstein classes (Corollaire \ref{R3}).

\end{tabular}
\end{center}

\newpage
\section{Introduction}

Les classes d'Eisenstein d'un sch\'ema ab\'elien (cf. \cite[Partie 5]{b} pour une d\'efinition)
ont une origine motivique d'apr\`es Kings (cf. \cite{ki1}). Dans le cas elliptique, celles-ci ont \'et\'e
intensivement \'etudi\'ees et on en conna\^it des propri\'et\'es remarquables, mais en dimension sup\'erieure,
peu de choses sont connues \`a leur sujet.
Par exemple, la non nullit\'e de ces classes constitue un probl\`eme ouvert.
Levin avait conjectur\'e que les courants qu'il construit dans $\cite{l}$ permettent de d\'ecrire le polylogarithme (cf. \cite[Partie 4.1]{b} pour une d\'efinition) et par suite les classes d'Eisenstein d'un sch\'ema ab\'elien. Cette conjecture \'etant  d\'emontr\'ee (cf. \cite{b}), on peut envisager d'utiliser cet outil pour
aborder l'\'etude des classes d'Eisenstein en dimension sup\'erieure. \\

Dans cet article, nous consid\'erons le cas particulier des familles de vari\'et\'es ab\'eliennes de Hilbert-Blumenthal
et nous utilisons le-dit outil pour exhiber des propri\'et\'es remarquables des classes d'Eisenstein. Le lecteur pourra consulter la note \cite{bn} pour un expos\'e concis des travaux pr\'esent\'es ici. Toutefois, dans la Proposition 3.1 et dans le Théorème 4.1 a) de \cite{bn}, les constantes rationnelles sont erronées. Les valeurs correctes de celles-ci sont données respectivement dans la Proposition 4.3 et le Théorème  5.2. On commence par donner une expression en coordonn\'ees de ces classes, 
au niveau topologique (cf. Proposition \ref{explicit_Eis}).
Celles-ci vivent dans un groupe de cohomologie de la base du sch\'ema ab\'elien, base qui est une vari\'et\'e de
Hilbert-Blumenthal. On \'etudie ensuite leur d\'eg\'en\'erescence en une pointe de la compactification de Baily-Borel de la base et l'on d\'emontre que les r\'esidus en cette pointe de ces classes s'expriment en termes de valeurs sp\'eciales d'une certaine
fonction $L$ attach\'ee au corps de nombres totalement r\'eel sous-jacent (cf. Th\'eor\`eme \ref{R1}).
On r\'eduit ce calcul \`a un calcul au niveau
topologique (cf. Proposition \ref{prig}) en utilisant de fa\c{c}on essentielle le Th\'eor\`eme de Burgos-Wildeshaus \cite[Theorem 2.9]{bw}. Ce r\'esidu \'etant un nombre rationnel, on en d\'eduit une preuve g\'eom\'etrique du Th\'eor\`eme de
Klingen-Siegel (cf. Corollaire \ref{R2}). \\

Plusieurs preuves de ce Th\'eor\`eme ont \'et\'e donn\'ees depuis 1962, date de la publication de la premi\`ere d\'emonstration dans l'article \cite{kl}.
Si notre approche pr\'esente une certaine analogie avec la preuve originale,
elle met toutefois en lumi\`ere un lien g\'eom\'etrique nouveau avec la famille modulaire de Hilbert-Blumenthal, via le polylogarithme. On remarque que les deux d\'emonstrations les plus r\'ecentes de ce r\'esultat dues \`a Nori \cite{no}
et Sczech \cite{sc} utilisent toutes deux la rationalit\'e d'une classe de cohomologie pour en d\'eduire la rationalit\'e d'une valeur sp\'eciale de la fonction $L$ en question. Notre approche pr\'esente donc \'egalement une certaine analogie avec les leurs.\\

On d\'eduit \'egalement de la forme particuli\`ere de la d\'eg\'erescence des classes d'Eisenstein un r\'esultat de non
nullit\'e pour certaines classes d'Eisenstein dans cette situation g\'eom\'etrique (cf. Corollaire \ref{R3}). La preuve utilise une \'equation fonctionnelle due \`a Siegel pour la fonction $L$ qui appara\^it dans l'expression du r\'esidu. On mentionne enfin que Kings dans la pr\'epublication \cite{ki2} \'etudie \'egalement la d\'eg\'en\'erescence des classes d'Eisenstein dans le cas Hilbert-Bimenthal.\\

On passe maintenant au contenu des diff\'erentes parties de cet article. \\
\begin{itemize}
\item[$\bullet$] Partie 2: On d\'efinit le sch\'ema ab\'elien que l'on consid\`ere \`a l'aide du formalisme de Pink (cf. \cite{p}). Cette approche nous est utile pour la suite. On d\'efinit \'egalement des sections de torsion et on esquisse
la construction de la compactification de Baily-Borel d'une composante connexe de la base.\\

\item[$\bullet$] Partie 3: Cette partie est consacr\'ee \`a la d\'efinition et \`a l'\'etude du morphisme r\'esidu en une pointe de la compactification de Baily-Borel de la base du sch\'ema ab\'elien pr\'ec\'edemment d\'efini. On y montre
en particulier que le but de ce morphisme est soit nul, soit canoniquement isomorphe \`a $\Q$. C'est ici qu'intervient
le Th\'eor\`eme de Burgos-Wildeshaus \cite[Theorem 2.9]{bw}. Son \'enonc\'e faisant appel au formalisme de Pink, la d\'efinition du sch\'ema ab\'elien dans ce langage nous est pr\'ecieuse. On \'etablit \'egalement un r\'esultat de rigidit\'e (cf. Proposition \ref{prig}) qui nous permet de r\'eduire le calcul du r\'esidu \`a un calcul au niveau topologique.\\

\item[$\bullet$] Partie 4: Il s'agit ici de calculer les diff\'erents objets qui interviennent dans la d\'efinition des courants de Levin (cf. \cite{l}). Le lecteur trouvera une explication des diff\'erents objets en jeu en d\'ebut de section. Ces calculs effectu\'es, on peut alors donner une formule
pour les classes d'Eisenstein, au niveau topologique, dans notre contexte (cf. Proposition \ref{explicit_Eis}). \\

\item[$\bullet$] Partie 5: Disposant d'une expression explicite au niveau topologique des classes d'Eisenstein et d'un r\'esultat de rigidit\'e qui nous permet d'effectuer le calcul au niveau topologique, on  d\'etermine, dans cette partie, le r\'esidu de ces classes en int\'egrant celles-ci le long du bord d'un voisinage ouvert de la pointe consid\'er\'ee. On
\'etablit alors le lien entre ces r\'esidus et des valeurs sp\'eciales d'une certaine fonction $L$
(cf. Notation \ref{fl} et Th\'eor\`eme \ref{R1}) et on en d\'eduit le Th\'eor\`eme de Klingen-Siegel (cf. Corollaire \ref{R2}) et un r\'esultat de non annulation (cf. Corollaire \ref{R3}).\\

\end{itemize}

\noindent {\bf Remerciements} \\

La d\'eg\'en\'erescence des classes d'Eisenstein des familles modulaires de Hilbert-Blumenthal
en des valeurs sp\'eciales de fonctions $L$ du corps du nombres totalement r\'eel sous jacent
avait \'et\'e conjectur\'ee par  J\"org Wildeshaus.
Je tiens \`a lui exprimer toute ma gratitude pour m'avoir fait part de cette intuition,
ainsi que
pour les entretiens qu'il a bien voulu m'accorder.
Je suis \'egalement ravi de remercier Guido Kings pour les discussions que nous avons partag\'ees sur
 les r\'esultats de cet article,
 lors de mon s\'ejour \`a l'Universit\'e de Regensburg. Enfin, j'adresse un grand merci à Daniel Barsky qui a eu 
 la gentillesse de bien vouloir répondre à mes questions sur les fonctions $L$ considérées dans ce travail.
 \\

\noindent{\bf Notations:} \\

Pour $X$ un sch\'ema de type fini, s\'epar\'e et lisse sur $\C$, on note: \\

\begin{tabular}{lp{13cm}}
$VSHM(X)\;$   & la cat\'egorie des $\Q$-variations de structures de Hodge mixtes polarisables et admissibles
                (cf. \cite{ka}) 
                 sur $X$, \\
$\overline{\mathbb{V}}$ & le syst\`eme local sous-jacent
          \`a $\mathbb{V}$ pour $\mathbb{V} \in Ob(VSHM(X))$,\\
$MHM(X)$ & la cat\'egorie des $\Q$-modules de Hodge alg\'ebriques mixtes sur $X$ (cf. \cite{s}).\\
         & \\
\end{tabular}

En particulier, $MHM(\text{Spec}(\C))$ est
la cat\'egorie des $\Q$-structures de Hodge mixtes  polarisables, cat\'egorie que l'on note simplement $SHM$.\\

On a un foncteur canonique
$VSHM(X) \to MHM(X)$ qui est exact, pleinement fid\`ele et gr\^ace auquel on identifie
$VSHM(X)$ \`a une sous-cat\'egorie pleine de $MHM(X)$.\\
Par construction de $MHM(X)$,
on dispose d'un foncteur $\For := real \circ rat \colon D^b MHM(X) \to D^b_c(X)$,
o\`u $rat$ est le foncteur d\'efini par M. Saito,
$real$ est le foncteur de Beilinson et
$D^b_c(X)$ la sous cat\'egorie pleine de la cat\'egorie
d\'eriv\'ee born\'ee des faisceaux de $\Q$-vectoriels sur $X(\C)$, \'equip\'e de la topologie transcendante,
dont les objets sont les complexes \`a cohomologie alg\'ebriquement constructible.\\

\begin{conv} $-$ 
Si $\mathbb{V} \in Ob(VSHM(X))$, $For(\mathbb{V})$ est le syst\`eme local $\overline{\mathbb{V}}$ convenablement
d\'ecal\'e. Dans cet article, {\it on ne tient pas compte du d\'ecalage}. \\
\end{conv}

On fixe pour la suite: \\

\begin{tabular}{cp{14cm}}
$L$ & un corps de nombres totalement r\'eel $L$ de dimension $g$, \\
$(\sigma_k)_{1 \leq k \leq g}$  & une \'enum\'eration des plongements r\'eels de $L$,\\
$\mathfrak{a}$ & un id\'eal entier de $L$,\\
$N$ & un nombre entier plus grand que 3.\\ \\
\end{tabular}

Soient $\mathfrak{b}$ un id\'eal fractionnaire de $L$ et $\mathfrak{p}$ un id\'eal premier de $L$. On note \\

\begin{tabular}{ll}
$\OO_L$ & l'anneau des entiers de $L$, \\
$\Tr_L$ & la trace de $L$, \\
 $\text{N}_L$ & la norme de $L$, \\
$d_L$ & le discriminant de $L$, \\
$\mathfrak{d}_L$ & la diff\'erente de $L$, \\
  $\mathcal{N}_L( \mathfrak{b} )$ & la norme de $\mathfrak{b}$, \\
 $\mathfrak{b}^{\vee}$ & le dual de $\mathfrak{b}$ relativement \`a $\Tr_L$, \\
$L_{\mathfrak{p}}$ & la compl\'etion de $L$ en la valuation associ\'ee \`a $\mathfrak{p}$, \\
$\OO_{L_{\mathfrak{p}}}$ & l'anneau des entiers de $L_{\mathfrak{p}}$, \\
$\mathfrak{b}_{\mathfrak{p}}$ & la compl\'etion de $\mathfrak{b}$ en la valuation associ\'ee \`a $\mathfrak{p}$.

\end{tabular}

\section{Donn\'ees g\'eom\'etriques}

On commence par d\'efinir le sch\'ema ab\'elien
que l'on consid\`ere dans cet article.
On utilise pour cela le formalisme et les r\'esultats de la Th\`ese de Pink \cite{p}.
On rappelle ensuite rapidement la construction de la compactification de Baily-Borel
de la base du sch\'ema ab\'elien pr\'ec\'edemment d\'efini, qui est une vari\'et\'e de
Hilbert-Blumenthal. On d\'efinit enfin des sections de torsion.

\begin{conv} $-$ 
Soit $H$ un groupe algébrique sur $L$ et soit $K \in \{\R,\C\}$.
Pour toute $K$-algèbre $A$, on dispose d'un isomorphisme de $K$-algèbres
           $ A \otimes_{\Q} L  \isom   A^g, \;
                                 a \otimes l  \mapsto  ( \sigma_k(l)  a )_{1 \leq k \leq g}$ 
au moyen duquel on identifie $(\Res_{L/\Q} H) \times_{\Q} K$ et $(H \times_{\Q} K )^g$.\\
\end{conv}	   

\subsection{D\'efinition de la base du sch\'ema ab\'elien}

Soit $i: \G_{m,\Q}  \hookrightarrow   \Res_{L/\Q}\G_{m,L}$ le morphisme qui à une $\Q$-algèbre $A$ fait correspondre le morphisme de groupes
$ A^\times \to (A \otimes L)^\times, \;a \mapsto a \otimes 1.$ On d\'efinit le $\Q$-sch\'ema en 
groupes r\'eductif $G$ par le diagramme cart\'esien suivant:

  $$\xymatrix{  \ar @{} [dr] |{\Box} G \; \ar@{^{(}->}[r] \ar[d] & \Res_{L/\Q} GL_{2,L} \ar[d]^{\Res_{L/\Q}(\mbox{d\'et})}  &\\ 
                \G_{m,\Q} \ar@{^{(}->}[r]_-i &  \Res_{L/\Q}\G_{m,L}\;. }$$

On pose $\h^{\pm} := \C-\R$ et on fait de $(\h^{\pm})^g$ un $G(\R)$-espace homogène \`a
gauche en faisant agir le groupe $G(\R) \subseteq GL_2(\R)^g$ à gauche sur $(\h^{\pm})^g$
via les homographies. On d\'efinit alors un morphisme 
$h : (\h^{\pm})^g \to \Hom(\Sp_{\C},G_{\C}) \subseteq \Hom(\Sp_{\C},(GL_{2,\C})^g )$ par:
$$\tau   \in (\h^{\pm})^g \mapsto 
\left[(z_1,z_2)  \in \Sp(\C) \mapsto 
   \left( \frac{i}{2 \;\text{Im}(\tau_k)} 
    \left( \begin{array}{cc} \overline{\tau_k} z_1 - \tau_k z_2 & -|\tau_k|^2 (z_1 - z_2) \\ 
					      z_1 - z_2 & -\tau_k z_1 + \overline{\tau_k}z_2 \end{array} \right)
      \right)_{1 \leq k \leq g}  \in  GL_{2}(\C)^g   \right].$$    
Alors, $(G, (\h^{\pm})^g)$  est une donnée de Shimura pure. Pour tout $\mathfrak{p}$ ideal premier de $L$, 
on pose 
$$ H(\mathfrak{a},N,\mathfrak{p}) :=
   \left\{ \left( \begin{array}{cc} a & b \\ c & d \end{array} \right) \in GL_2(L_{\mathfrak{p}}) \; 
   \left| \; \begin{array}{l}
              a,d \in 1 + N \OO_{L_{\mathfrak{p}}}, \\
	      c \in N (\mathfrak{d}_L\mathfrak{a}^2)_{\mathfrak{p}}, \;
              b \in N (\mathfrak{d}_L^{-1}\mathfrak{a}^{-2})_{\mathfrak{p}}, \\
	      ad - bc \in \OO_{L_{\mathfrak{p}}}^{\times} \end{array}
						         \right. \right\}.$$
Le groupe
$ H(\mathfrak{a},N) := \left( \prod\limits_{\mathfrak{p}} \; H(\mathfrak{a},N,\mathfrak{p}) \right) \cap G(\mathbb{A}_{\Q,f})$
est un sous-groupe compact ouvert net de $G(\mathbb{A}_{\Q,f})$.\\ 

On note $S$ la vari\'et\'e alg\'ebrique complexe quasi-projective lisse
$M^{H(\mathfrak{a},N) } (G, (\h^{\pm})^g )$ (cf. \cite[p. 103]{g} pour une interpr\'etation modulaire de $S$).
 L'inclusion canonique
du produit de $g$ copies du demi-plan  de Poincar\'e sup\'erieur $\h$ dans
$(\h^{\pm})^g$ induit une immersion ouverte
$$ \Gamma(\mathfrak{a},N) \backslash  \h^g \hookrightarrow S^{an}  ,$$
o\`u $\Gamma(\mathfrak{a},N) :=
      \left\{ \left( \begin{array}{cc} a & b \\ c & d \end{array} \right) \in SL_2(L) \; \left| \;
                  a,d \in 1 + N \OO_{L}, \; c \in N \mathfrak{d}_L\mathfrak{a}^2,
		                                      \; b \in N \mathfrak{d}_L^{-1}\mathfrak{a}^{-2} 
						      \right. \right\} \subseteq G(\Q),$
qui identifie $ \Gamma(\mathfrak{a},N) \backslash  \h^g$
 \`a une composante connexe de $S^{an}$. 
  La vari\'et\'e analytique complexe
$ \Gamma(\mathfrak{a},N) \backslash  \h^g$ est la vari\'et\'e analytique complexe associ\'ee
\`a une vari\'et\'e alg\'ebrique complexe quasi-projective lisse canonique not\'ee $S^0$. 
Par d\'efinition de la structure de vari\'et\'e alg\'ebrique de $S$, 
$S^0$ est une composante connexe de $S$.

\subsection{D\'efinition de la famille de vari\'et\'es ab\'eliennes}

Soit $V:= \Res_{L/ \Q}(\G_{a,L} \oplus \G_{a,L})$.
On pose $P:= V \rtimes G$ le produit semi-direct associ\'e \`a la restriction \`a $G$ de l'action
standard de $\Res_{L/ \Q}GL_{2,L}$ sur $V$.

\begin{lemme}{\emph{(comparer à [Wi1-V-Lemma 1.1])}} $-$
Soient $\tau \in (\h^{\pm})^g$. Le morphisme $h(\tau)$
définit une $\Q$-structure de Hodge pure de poids $-1$ sur $V(\Q)$.
Sous l'identification $V_{\C} = ( \G_{a,\C} \oplus \G_{a,\C} )^g$ la filtration de Hodge est donnée par:
$$ \begin{array}{l}
    F^1(V_{\C}) = 0, \\
    F^0(V_{\C}) = \left< \left( \begin{array}{c} \tau_1 \\ 1 \end{array} \right) \right>_{\C} \times  \dots \times
 		 \left< \left( \begin{array}{c} \tau_g \\ 1 \end{array} \right) \right>_{\C}, \\
     F^{-1}(V_{\C}) = V_{\C}.
\end{array} $$
\end{lemme}

\begin{dem} $H^{0,-1}(V_{\C})$ est le sous-espace propre associé au caractère $(z_1,z_2) \mapsto z_2$.
Un calcul permet alors de conclure. \end{dem}

On définit alors la donnée de Shimura $(P,\X)$ comme \'etant l'extension unipotente de $(G, (\h^{\pm})^g )$ par $V$. Comme $\text{Lie}(V)$ est de type $\{(-1,0) , (0,-1) \}$, on a la description suivante de $\X$:
$$ \X = \left\{ (\psi,\tau) \in \Hom(\Sp,P_{\R}) \times (\h^{\pm})^g \; \left| \;
    h(\tau) = p \circ \psi  \right. \right\}, $$
où $p : P \to G$ est la projection canonique.
\`A l'aide des isomorphismes donnés par [Wi1-V-Lemma 1.2] et [Wi1-V-Corollary 1.4.b], on construit
un biholomorphisme:
$$ \C^g \times (\h^{\pm})^g \isom \X .$$

On introduit le sous-groupe ouvert compact de $P(\mathbb{A}_{\Q,f})$
$$ K(\mathfrak{a},N) :=
   \left( \prod\limits_{\mathfrak{p}} \; (\mathfrak{a}^{\vee})_{\mathfrak{p}} \; \oplus
          \prod\limits_{\mathfrak{p}} \; \mathfrak{a}_{\mathfrak{p}} \right)
      \rtimes H(\mathfrak{a},N).$$

On note $A$ la vari\'et\'e alg\'ebrique complexe quasi-projective lisse
$M^{K(\mathfrak{a},N) } (P, \X )$. La projection canonique $ p:P \to G$ 
induit un morphisme $A \to S$, not\'e (abusivement) \'egalement $p$, 
qui est le morphisme structural d'un sch\'ema ab\'elien de dimension relative pure $g$. \\

La restriction de $p^{an} \colon A^{an} \to S^{an}$ au dessus $(S^0)^{an} =  \Gamma(\mathfrak{a},N) \backslash  \h^g$ est donn\'ee par
la projection canonique
$$ q \colon A^0 := \Lambda(\mathfrak{a},N)  \backslash \X^{+} \to (S^0)^{an},$$
o\`u $\X^{+} := \C^g \times \h^g$ et
$\Lambda(\mathfrak{a},N) := \left(  \mathfrak{a}^{\vee} \oplus \mathfrak{a} \right)
                             \rtimes \Gamma(\mathfrak{a},N)$
agit sur $ \X^{+}     $ par l'action qui à
$$\left(  \left( (a',a) , \left( \begin{array}{cc}
                                   \alpha & \beta \\
                                   \gamma & \delta
                        \end{array} \right) \right) ,
     ( z, \tau ) \right) \in
    \Lambda(\mathfrak{a},N) \times ( \C^g \times  \h^g)$$ fait correspondre:
$$\left(
      \frac{z_k}{\sigma_k(\gamma)  \tau_k + \sigma_k(\delta) } + \sigma_k(a') -
                 \sigma_k(a) \left(
		  \frac{ \sigma_k(\alpha) \tau_k + \sigma_k(\beta)}{\sigma_k(\gamma) \tau_k + \sigma_k(\delta)}
                \right)  , \frac{ \sigma_k(\alpha) \tau_k + \sigma_k(\beta)}{\sigma_k(\gamma) \tau_k +
		\sigma_k(\delta)} \right)_{1 \leq k \leq g} \in \C^g \times  \h^g.$$
On a donc l'identification suivante:  $\left( p_{|S^0} \colon A_{|S^0} \to S^0 \right)^{an} = (q \colon A^0 \to (S^0)^{an})$.

\subsection{La pointe $\infty$}

On rappelle succinctement comment est construite la compactification de Baily-Borel de $S^0$.
On d\'efinit une
pointe de $\Gamma(\mathfrak{a},N)$ comme \'etant 
une orbite d'un point de $\mathbb{P}^1(L)$ sous l'action 
de $\Gamma(\mathfrak{a},N) \subseteq SL_2(L)$ induite par les homographies.
On adjoint à $S^0$ l'ensemble des pointes de  $\Gamma(\mathfrak{a},N)$ pour obtenir un ensemble que l'on munit de la topologie de Satake. Sur cet espace topologique, on dispose
d'une structure d'espace analytique complexe normal compact (non lisse si $L \ne \Q$) compatible avec
celle de $(S^0)^{an}$. 
\`A l'aide des s\'eries de Poincar\'e, 
on construit un plongement de cet espace analytique dans un espace projectif.
Il h\'erite ainsi d'une structure de vari\'et\'e alg\'ebrique complexe normale et projective. On note  $(S^0)^*$ la vari\'et\'e alg\'ebrique obtenue.
 Par d\'efinition m\^eme de la structure alg\'ebrique de $S^0$, $S^0$ est un ouvert dense de $(S^0)^*$ et
 la vari\'et\'e $(S^0)^*$ est appel\'ee compactification de Baily-Borel de $S^0$.  \\

On d\'ecrit maintenant
un système fondamental de voisinages dans $(S^0)^*$ (pour la structure analytique)
de la pointe $\infty$, pointe qui correspond à l'orbite de $[1:0] \in \mathbb{P}^1(L)$ sous l'action de $\Gamma(\mathfrak{a},N)$.
 Le stabilisateur de l'$\infty$ dans $\Gamma(\mathfrak{a},N)$ est
$$ \Gamma(\mathfrak{a},N,\infty)  := \left\{ \left(
   \begin{array}{cc}
   \varepsilon & a \\
    0 & \varepsilon^{-1}
   \end{array} \right) \in \Gamma(\mathfrak{a},N) \; \left| \; \varepsilon \in U_{L,N} \mbox{ et }
        a \in N \mathfrak{d}_L^{-1}\mathfrak{a}^{-2} \right.
     \right\},
$$
o\`u $U_{L,N}$ est le sous-groupe des unit\'es de $\OO_L$ congrues \`a 1 modulo $N$.
On v\'erifie que $ \Gamma(\mathfrak{a},N,\infty)$ agit sur
   $$ V_r : =  \left\{ \tau  \in \h^g \; \left| \; \prod\limits_{i=1}^g  \text{Im}(\tau_j) > r 
   \right. \right\},$$
pour tout $r \in \R^{>0}.$ Alors pour $r \gg 0$, le morphisme canonique
$ \Gamma(\mathfrak{a},N,\infty) \backslash V_r \to S^0$ est une immersion ouverte et l'ensemble
$ \left\{ (\Gamma(\mathfrak{a},N,\infty) \backslash V_r) \cup \{ \infty \} \; \left| \; r \gg 0 \right. \right\}$
forme une base de voisinages ouverts de l'$\infty$ dans $S^0$.

\subsection{Sections de torsion}

Soient $a' \in N^{-1} \mathfrak{a}^{\vee}$, $a \in N^{-1} \mathfrak{a}$. Le morphisme analytique
$$\begin{array}{ccccc}
   x_{a',a} & \colon & (S^0)^{an} = \Gamma(\mathfrak{a},N) \backslash \h^g & \to &
   \Lambda(\mathfrak{a},N) \backslash ( \C^g \times  \h^g) \simeq (A_{|S^0})^{an} \\
   & & \left[ \tau \right] & \mapsto & \left[ ( \sigma_k(a') + \sigma_k(a) \tau_k)_{1 \leq k \leq g} , \tau \right].
\end{array}$$
est le morphisme analytique associ\'e \`a un unique morphisme alg\'ebrique $S^0 \to A_{|S^0}$ qui est une
section de $N$-torsion de $p_{|S^0} \colon A_{|S^0} \to S^0$ que l'on note (abusivement) \'egalement $x_{a',a}$.\\

\section{Le morphisme r\'esidu en l'$\infty$}

Cette article traite de la d\'eg\'en\'erescence des classes d'Eisenstein (cf. \cite[D\'efinition 5.3]{b}) du sch\'ema ab\'elien
$p_{|S^0} \colon A_{|S^0} \to S^0 $ en la pointe $\infty$, i.e. de la d\'etermination du r\'esidu
en l'$\infty$ des classes d'Eisenstein. \\

Dans cette partie, on d\'efinit le morphisme r\'esidu
en l'$\infty$, de source un groupe d'extensions de $MHM(S^0)$ et de but un groupe d'extensions de $SHM$ qui est soit
trivial, soit canoniquement isomorphe \`a $\Q$. On explique enfin que le calcul du r\'esidu en l'$\infty$ se r\'eduit
\`a un calcul au niveau topologique. \\

On note: \\

\begin{tabular}{cp{13cm}}
$j : S^0 \hookrightarrow (S^0)^*$ & l'immersion ouverte de $S^0$ dans $(S^0)^*$, \\
$i_{\infty} : \infty  \hookrightarrow (S^0)^*$ & l'immersion fermée de la pointe $\infty$
             dans $(S^0)^*$,\\
$\hh$ & $:= (R^1 (p_{|S^0})_* \Q)^{\vee} \in Ob(VSHM(S^0))$,\\
$ l$  & un entier positif.
\end{tabular}

\subsection{D\'efinition du morphisme r\'esidu en l'$\infty$}

On d\'efinit le morphisme $\Res_{\infty}^l$ par le diagramme commutatif suivant.
$$\!\!\!\xymatrix{
\Ext^{2g-1}_{MHM(S^0)}( \Q(0), (\Sym^l \hh)(g))
 \ar@{=}[r] \ar[ddd]_-{\Res_{\infty}^l}    &
 \Hom_{D^b MHM(S^0)}( \Q(0), (\Sym^l \hh)(g)[2g-1])
 \ar@{=}[d]^-{\; \text{(adjonction)}}  \\
 &
 \Hom_{D^b MHM((S^0)^*)}( \Q(0), j_*(\Sym^l \hh)(g)[2g-1])
 \ar[d]^-{\; 1 \to (i_{\infty})_* (i_{\infty})^*}\\
 &
  \Hom_{D^b MHM((S^0)^*)}( \Q(0), (i_{\infty})_* (i_{\infty})^* j_*(\Sym^l \hh)(g)[2g-1])
   \ar@{=}[d]^-{\; \text{(adjonction)}}  \\
  \Hom_{SHM} ( \Q(0) , H^{2g-1} (i_{\infty})^* j_*(\Sym^l \hh)(g))
  & \Hom_{D^b SHM}( \Q(0), (i_{\infty})^* j_*(\Sym^l \hh)(g)[2g-1])
  \ar[l]^-{\; H^0}
} \\
$$
$\;$ \\

\subsection{Le but du morphisme r\'esidu en l'$\infty$}

On d\'etermine \`a pr\'esent la $\Q$-structure de Hodge  $H^{2g-1} (i_{\infty})^* j_*(\Sym^l \hh)(g)$. \\

Par construction de $S =M^{H(\mathfrak{a},N) } (G, (\h^{\pm})^g )$, on dispose d'un foncteur canonique
qui associe \`a une repr\'esentation alg\'ebrique rationnelle de $G$ une $\Q$-variation de structures de Hodge
mixtes polarisable sur $M^{H(\mathfrak{a},N) } (G,(\h^{\pm})^g )(\C)$. De plus, d'apr\`es \cite[II - Theorem 2.2]{w},
les variations issues de cette construction sont admissibles. On en d\'eduit un foncteur tensoriel canonique
$$ \mu_{H(\mathfrak{a},N)} \colon \text{Rep}_{\Q} G \to MHM(S),$$
o\`u $\text{Rep}_{\Q} G$ est la cat\'egorie des repr\'esentations alg\'ebriques rationnelles de $G$· \\

D'apr\`es la d\'efinition de $G$, on a un morphisme canonique $G \to \G_{m,\Q}$, not\'e $\chi$, qui fournit
une action de $G$ sur $\G_{a,\Q}$. On v\'erifie que $\mu_{K(\mathfrak{a},N)}(\chi) = \Q(1)$.
Le groupe $G$ agit sur $V:= \Res_{L/ \Q}(\G_{a,L} \oplus \G_{a,L})$ par restriction de l'action standard
de $\Res_{L/ \Q} GL_{2,L} $. L'image par $\mu_{H(\mathfrak{a},N)}$ de cette repr\'esentation
est $ (R^1 p_{*} \Q)^{\vee}$ ($p$ est le morphisme structural du sch\'ema ab\'elien $A \to S$).
Ainsi $(\Sym^l (R^1 p_{*} \Q)^{\vee})(g)$ est dans l'image du foncteur $\mu_{H(\mathfrak{a},N)}$.\\

Burgos et Wildeshaus ont
\'etabli un Th\'eor\`eme \cite[Theorem 2.9]{bw} qui permet de calculer
la $\Q$-structure de Hodge  $ H^{2g-1} i_{\infty}^* j_* \; \mu_{H(\mathfrak{a},N) } (W),$
pour $W \in \text{Rep}_{\Q} G$,
en termes de cohomologies de deux groupes:
l'un algébrique unipotent, l'autre arithmétique.
L'énoncé du résultat
précis requiert la description de la compactification de Baily-Borel dans le formalisme
de Pink. Pour cela, on renvoie aux chapitres 4 et 6 de la Thèse de Pink [P] ou encore,
pour un résumé, à [BW-p. 365 et 366]. On explicite les différents objets qui interviennent dans la description
de la pointe $\infty$ de $B$ en suivant les notations de [BW- p. 365 -- 367]
 et aussi les groupes qui interviennent dans ce cas dans le Théorème [BW-Theorem 2.9]. \\

\begin{itemize}
\item[$\bullet$] 
           Soit $Q$ le sous-groupe parabolique admissible de $G$ qui est le produit fibré du
           sous-groupe de Borel standard de $\Res_{L/\Q} GL_{2,L}$
           $$ \left\{ \left( \begin{array}{cc} * & * \\ 0 & * \end{array} \right) \right\} $$
           et de $\mathbb{G}_{m,\Q}$.
\item[$\bullet$] Le sous-groupe normal canonique de $Q$, $P_1$, est le produit fibré du sous-groupe
            de $\Res_{L/\Q} GL_{2,L}$
            $$ \left\{ \left( \begin{array}{cc} * & * \\ 0 & 1 \end{array} \right) \right\} $$ 
            et de $\mathbb{G}_{m,\Q}$ qui s'identifie à
            $ \Res_{L/\Q} \mathbb{G}_{a,L} \rtimes \mathbb{G}_{m,\Q}. $
\item[$\bullet$] Le radical unipotent $W_1$ de $P_1$ (et de $Q$) est le produit fibré du sous-groupe de
             $\Res_{L/\Q} GL_{2,L}$
               $$ \left\{ \left( \begin{array}{cc} 1 & * \\ 0 & 1 \end{array} \right) \right\} $$ 
             et de $\mathbb{G}_{m,\Q}$, c'est à dire
             $ \Res_{L/\Q} \mathbb{G}_{a,L} .$
\item[$\bullet$]  On note $\pi \colon P_1 \to G_1 := P_1 / W_1 = \mathbb{G}_{m,\Q}$ la projection canonique.
\item[$\bullet$]   On fixe $K:=H(\mathfrak{a},N)$ sous-groupe compact ouvert net de $G(\mathbb{A}_{\Q,f})$,
                   $g:=1 \in G(\mathbb{A}_{\Q,f})$ et on pose
		    $K_1:= P_1( \mathbb{A}_{\Q,f}) \cap H(\mathfrak{a},N).$
\item[$\bullet$] La composante de bord dans la compactification de Baily-Borel 
            qui correspond à la pointe $\infty$ est la donnée de Shimura
            $ (\mathbb{G}_{m,\Q}, \text{Isom}(\Z,\Z(1)))$ (cf. \cite[Example 2.8]{p}).
\item[$\bullet$]  On a donc un morphisme
             $i_{\mathbb{G}_{m,\Q},H(\mathfrak{a},N),1} \colon M^{\pi(K_1)}(\mathbb{G}_{m,\Q}
	      ,\text{Isom}(\Z,\Z(1))) \to 
	     (S^0)^*.$
\item[$\bullet$] On calcule 
                 $H_C := \text{Cent}_{Q(\Q)}(\{ -1,1 \}) \cap W_1(\mathbb{A}_{\Q,f}) \cdot H(\mathfrak{a},N)$
		 et $ \overline{H_C} := \pi(H_C)$ pour trouver
$$ H_C = \left\{ \left( \begin{array}{cc}
              \varepsilon & * \\ 0 & \varepsilon^{-1} \end{array} \right) \; \left| \;
             \varepsilon \in U_{L,N}  \right. \right\} \quad \text{ et } \quad
             \overline{H_C} =
            \left\{ \left[ \left( \begin{array}{cc}
              \varepsilon & * \\ 0 & \varepsilon^{-1} \end{array} \right) \right] \;  \left| \;
             \varepsilon \in U_{L,N} \right. \right\} \simeq U_{L,N}.$$
	    
\item[$\bullet$] 
             Enfin,    
             $\Delta \backslash  M^{\pi(K_1)}(\mathbb{G}_{m,\Q}  , \text{Isom}(\Z,\Z(1)))$
	      (cf. \cite[p. 367]{bw} pour la d\'efinition de $\Delta$) est réduit à un point
             et l'immersion  induite par $i_{\mathbb{G}_{m,\Q},H(\mathfrak{a},N),1} 
	     \colon M^{\pi(K_1)}(\mathbb{G}_{m,\Q} ,\text{Isom}(\Z,\Z(1))) \to 
	     (S^0)^*$
                $$ \Delta \backslash  M^{\pi(K_1)}(\mathbb{G}_{m,\Q}, \text{Isom}(\Z,\Z(1)))
		\hookrightarrow   (S^0)^* $$
             correspond à
             $ i_{\infty} :  \infty  \hookrightarrow (S^0)^* .$

\end{itemize}
\vspace{0.3cm}
On peut alors énoncer le résultat suivant.

\begin{thm}{(cas particulier de [BW-Theorem 2.9])} $-$
$$ H^{2g-1} i_{\infty}^*j_* (\Sym^l \hh) (g) =
   \underset{p+q=2g-1}{\oplus} \; \mu_{\pi(K_1)}
    \circ H^p(\overline{H_C},H^q(W_1,Res^G_Q ( (\Sym^l V) \otimes \chi^g))) .$$
\end{thm}

La dimension cohomologique du groupe $W_1$ (resp. du groupe abélien libre de rang $g-1$ sans torsion
 $\overline{H_C}$) est $g$ (resp. $g-1$). Ainsi, on a:
$$ H^{2g-1}i_{\infty}^*j_* (\Sym^l \hh) (d) =
    \mu_{\pi(K_1)}
    \circ H^{2g-1}(\overline{H_C},H^{g}(W_1,Res^G_Q ( (\Sym^l V) \otimes \chi^g))). $$

\begin{prop}{} $-$ \label{deg}
$$
    H^{2g-1}i_{\infty}^*j_* (\Sym^l \hh)(g)    =
    \left\{ \begin{array}{cl} \Q(0) & \mbox{ si } g \mbox{ divise } l, \\
                               0  & \mbox{ sinon }.
            \end{array} \right.
$$
\end{prop}

\newpage

\begin{dem}\\

\begin{itemize}
\item[a) ] Calcul de $H^{g}(W_1,Res^G_Q ( \Sym^l V \otimes \chi^g))$ \\

On commence par remarquer que $W_1$ n'agit pas sur $\chi^g$.
D'après [Kn-Thm 6.10], on a un isomorphisme $(Q/W_1)$-équivariant:
$$  H^{g}(W_1,\Res^G_Q ( \Sym^l ) ) \simeq H_0(W_1,\Res^G_Q ( \Sym^l V)) \otimes
           \overset{g}{\wedge} (\text{Lie} W_1)^{\vee}. $$
On étend ensuite les scalaires à $\overline{\Q}$, une clôture algébrique de $\Q$. On a l'isomorphisme canonique suivant:
  $$ H_0(W_1,\Res^G_Q ( \Sym^l V )) \otimes_{\Q} \overline{\Q} =
      H_0(W_1( \overline{\Q}), \Sym^l V(\overline{\Q}) ) .$$
On note que:
$$ \begin{array}{ccl}
    W_1( \overline{\Q}) & = & \displaystyle \prod\limits_{i=1}^g
                             \left\{ \left( \begin{array}{cc} 1 & \alpha \\
                         0 & 1 \end{array} \right) \; \left| \; \alpha \in \overline{\Q} \right. \right\}\\
    V(\overline{\Q}) & = &   \displaystyle \prod\limits_{i=1}^g
    \left\{ a_i X_i + b_i Y_i  \;  \left| \; a_i,b_i \in \overline{\Q}  \right. \right\} ,
\end{array}$$
où $X_i = (1,0), \; Y_i = (0,1) \in \overline{\Q}^2$. On a ainsi une base canonique pour
$\Sym^l V(\overline{\Q})$:
$$ \left( X_1^{m_1} \dots X_g^{m_g}Y_1^{n_1} \dots Y_g^{n_g} \right)_{(m_1, \dots ,m_g,n_1, \dots ,n_g)} $$
indexée par les $2g$-uplets d'entiers positifs $(m_1, \dots ,m_g,n_1, \dots ,n_g)$ tels que:
$$ m_1 + \dots  + m_g + n_1 + \dots + n_g = l .$$
Soit $\alpha \in \overline{\Q}$ et soit
$w: =  \left( \begin{array}{cc} 1 & \alpha \\
                         0 & 1 \end{array} \right)$.
Alors, $w . X_i = X_i$ et $w Y_i = \alpha X_i + Y_i$ et on v\'erifie, \`a l'aide de cette remarque, que:
$$ H_0(W_1( \overline{\Q}), \Sym^l V(\overline{\Q}) ) =
   \Sym^l V(\overline{\Q}) /
\left<  \{ X_1^{m_1} \dots X_g^{m_g}Y_1^{n_1} \dots Y_g^{n_g} \; \left| \; \exists i \mbox{ tel que } m_i
 \right. \not= 0 \}
          \right>_{\overline{\Q}} .$$

\item[b) ] Calcul de $H^{g-1}(\overline{H_C},  ( \Sym^l V (\overline{\Q}) /
\left<  \{ X_1^{m_1} \dots X_g^{m_g}Y_1^{n_1} \dots Y_g^{n_g} \; \left| 
     \; \exists i \mbox{ tel que } m_i \not= 0 \right. \}
        \right>_{\overline{\Q}} ) \otimes \chi^g ).$ \\

Tout d'abord, $\overline{H_C}$ n'agit ni sur 
$\chi^g$ (les éléments de $\overline{H_C}$ sont de déterminant $1$) ni sur
$\overset{g}{\wedge} (\text{Lie} W_1)^{\vee}$ (le carré d'une unité est de norme 1).
Le groupe $\overline{H_C}$ étant isomorphe à $\Z^{g-1}$, en prenant une résolution
de Koszul, on obtient un isomorphisme entre:
$$
H^{g-1}(\overline{H_C},  ( \Sym^l V(\overline{\Q}) /
\left<  \{ X_1^{m_1} \dots X_g^{m_g}Y_1^{n_1}..Y_g^{n_g} \; \left| \; \exists i \mbox{ tel que }
      m_i \not= 0  \right. \}
        \right>_{\overline{\Q}} ))$$
et
$$ H_0 (\overline{H_C},  ( \Sym^l V(\overline{\Q}) /
\left<  \{ X_1^{m_1} \dots X_g^{m_g}Y_1^{n_1} \dots Y_g^{n_g} \; \left| \; \exists i \mbox{ tel que }
 m_i \not= 0 \} \right.
        \right>_{\overline{\Q}} ) .$$

Soient $(n_1, \dots ,n_g)$ un $g$-uplet d'entiers positifs tels que
$ n_1 + \dots + n_g = l $
et
$$ h = \left[ \left( \begin{array}{cc}
              \varepsilon & * \\ 0 & \varepsilon^{-1} \end{array} \right) \right]  \in \overline{H_C}$$
alors,
$ h . [ Y_1^{n_1} \dots Y_g^{n_g} ] =
      [ \sigma_1(\varepsilon^{-1})^{n_1}  \dots \sigma_g(\varepsilon^{-1})^{n_g}  Y_1^{n_1} \dots Y_g^{n_g} ].$
On remarque que si les $n_i$ ne sont pas tous égaux, alors il existe
$\varepsilon \in U_{L,N}$ tel que
 $ \sigma_1(\varepsilon^{-1})^{n_1}  \dots \sigma_g(\varepsilon^{-1})^{n_g} \not= 1. $
Ceci peut se voir en utilisant une $\Z$-base $(u_1, \dots ,u_{g-1})$ de $U_{L,N}$
et le résultat de théorie des nombres classique qui affirme que:
$$ \det \left( \begin{array}{ccc}
              \log |\sigma_1(u_1)|  & \dots & \log |\sigma_{g-1}(u_1)| \\
                        \vdots     & \dots &   \vdots  \\
               \log |\sigma_{1}(u_{g-1})| & \dots &    \log |\sigma_{g-1}(u_{g-1})|
               \end{array} \right) \not= 0 .$$

\item[c) ] Conclusion \\

On rassemble les résultats précédents. On définit $W$ comme suit. \\

\begin{tabular}{llp{11.2cm}}
$\bullet$ & Cas o\`u $g$ divise $l$: & On pose $\lambda := l /g \in \N $ et
on d\'efinit 
 $W$ comme \'etant le sous-espace de $\Sym^l V(\overline{\Q})$ engendré par les
$ X_1^{m_1} \dots X_g^{m_g}Y_1^{n_1} \dots Y_g^{n_g} $ tels que
$   (n_1, \dots ,n_g) \not= (\lambda, \dots ,\lambda) $.\\
& & \\
\end{tabular}\\

\begin{tabular}{llp{10cm}}
$\bullet$ & Cas o\`u $g$ ne divise pas $l$: & On pose $ W:= \Sym^l V(\overline{\Q})$.\\
& & \\
\end{tabular}

On a prouvé que:
$$ H^{2g-1 }i_{\infty}^*j_* (\Sym^l \hh)(d)  \otimes_{\Q} \overline{\Q}  =
     \Sym^l V (\overline{\Q}) / W \otimes \overset{g}{\wedge} (\text{Lie} W_1)^{\vee} \otimes \chi^g .$$
L'action de $G_1$ sur $(\text{Lie} W_1)^{\vee} \otimes \chi^g$ et sur $\Sym^l V (\overline{\Q}) / W$
est triviale. Enfin,
$$ \dim_{\overline{\Q}}  \;
 \Sym^l V(\overline{\Q}) / W =
 \left\{ \begin{array}{l} 1 \mbox{ si } g \mbox{ divise }  l ,\\
                         0 \mbox{ sinon. } \end{array} \right.$$
\end{itemize}
\end{dem}

\begin{rema} $-$ \label{remproj} Si $g$ divise $l$, alors le quotient $\Sym^l V(\overline{\Q}) / W$
(engendr\'e par la classe $[Y_1^{l/g} \dots  Y_g^{l/g} ]$) et la projection canonique
$\Sym^l V(\overline{\Q})  \to \Sym^l V(\overline{\Q}) / W$ sont
d\'efinis sur $\Q$.
\end{rema}

\subsection{Rigidit\'e de la d\'eg\'enerescence} \label{rig}

On suppose ici que $g$ divise $l$. La construction de $\Res_{\infty}^l$ admet un analogue \'evident
au niveau topologique (i.e. dans la th\'eorie $D^b_c( \cdot)$) qui permet de d\'efinir un morphisme
 $$ \overline{Res^l_{\infty}} \colon
H^{2g-1}_{\text{Betti}}( S^0(\C) , \overline{ (\Sym^l \hh) (g) } ) \to
H^{ 2g-1}_{\text{Betti}} (  \infty  , i_{\infty}^* Rj_* \overline{(\Sym^l \hh)  (g)} ).$$

Ce morphisme s'ins\`ere dans le diagramme suivant
$$
\xymatrix{
Ext^{2g-1}_{MHM(S^0)}( \Q(0)  , (\Sym^l \hh) (g))  \ar[d]_-{Res^l_{\infty}} \ar[r]^-{For} &
H^{2g-1}_{\text{Betti}}( S^0(\C) , \overline{\Sym^l \hh (g) } ) \ar[d]^-{\overline{Res^l_{\infty}}} \\
Hom_{SHM}(\Q(0),H^{ 2g-1} i_{\infty}^* j_* (\Sym^l \hh) (g) ) \ar[r]^-{For}
\ar@{=}[d]_-{\mbox{\tiny (cf  Proposition  \ref{deg})}} &
H^{ 2g-1}_{\text{Betti}} ( \infty , i_{\infty}^* Rj_* (\overline{\Sym^l \hh ) (g)} )
  \ar@{=}[d]^-{\mbox{\tiny (cf.  Proposition  \ref{deg})}} \\
Hom_{SHM}(\Q(0),\Q(0)) \ar@{=}[r]^-{} & \Q }$$
qui est commutatif. En effet, le formalisme des 6 foncteurs de $D^b MHM(\cdot)$ et celui de
$D^b_c( \cdot)$ sont compatibles via le foncteur $For$. On en d\'eduit la Proposition suivante.

\begin{prop} $-$ \label{prig}
Pour tout $c \in Ext^{2g-1}_{MHM(S^0)}( \Q(0) , (\Sym^l \hh) (g) )$,
$$Res^l_{\infty}(c) =  \overline{Res^l_{\infty}} \circ For (c) \in \Q.$$
\end{prop}

\section{Calcul des classes d'Eisenstein au niveau topologique}

Dans cette partie, on d\'etermine les classes d'Eisenstein (cf. \cite[Partie 5]{b} pour la d\'efinition)
 du sch\'ema ab\'elien $p_{|S^0} \colon A_{|S^0} \to S^0$ au niveau topologique en utilisant
 le r\'esultat principal de \cite{b} (cf. \cite[Corollaire 4.7 ]{b}).
Ce dernier assure que les courants d\'efinis par Levin dans \cite{l} permettent de d\'ecrire le polylogarithme
de $p_{|S^0} \colon A_{|S^0} \to S^0$ au niveau topologique, et par suite (cf. \cite[Partie 5]{b}) les classes d'Eisenstein, toujours au niveau topologique. \\

Il s'agit donc ici d'expliciter les courants de Levin dans notre situation g\'eom\'etrique.
Ces courants sont des s\'eries de formes diff\'erentielles not\'ees $g'_{a,\gamma}$ dans \cite{l}. Pour la
d\'efinition pr\'ecise de ces formes, on renvoie le lecteur \`a \cite{l}. On en donne ci-dessous une esquisse
tout en pr\'esentant le contenu de cette quatri\`eme partie.
\\

La construction de Levin vaut pour une famille de tores analytiques complexes munie d'une polarisation
(cf. \cite[1.1.2]{l}).
Dans la partie \ref{trivialisation}, on d\'efinit cette polarisation  et on explicite
 le pullback de
$$ (p_{|S^0} \colon A_{|S^0} \to (S^0))^{an} =
(q \colon \Lambda(\mathfrak{a},N) \backslash ( \C^g \times  \h^g) \to \Gamma(\mathfrak{a},N) \backslash \h^g)$$
par la projection canonique $ \h^g \to \Lambda(\mathfrak{a},N) \backslash \h^g$, not\'e 
$q' \colon A' \to \h^g$. Pour la famille de tores analytiques complexes $q' \colon A' \to \h^g$, on dispose de coordonn\'ees globales pour expliciter $g'_{a,\gamma}$. C'est la raison pour laquelle
on consid\`ere plutot cette famille que $(p_{|S^0} \colon A_{|S^0} \to (S^0))^{an}$. \\

On donne alors une trivialisation $\cinf$-r\'eelle $\Psi$ de la famille de tores r\'eels $q' \colon A' \to \h^g$
$$ \Psi \colon (\Pi_{\R} / \Pi) \times \h^g \to A',$$
o\`u $\Pi := \mathfrak{a}^{\vee} \oplus \mathfrak{a}$.
Le syst\`eme local constant $(R^1 q'_* \Z)^{\vee}$ s'identifie \`a $\Pi$ et la
polarisation (principale) est donn\'ee par un accouplement $< \cdot \; ,  \cdot > \colon \Pi \wedge \Pi \to 2 \pi i \Z$. En outre, $(R^1 q'_* \Z)^{\vee}$ est muni d'une structure de variation de structures de Hodge pure de poids $-1$.
On note $E$ le fibr\'e holomorphe sur $\h^g$ associ\'e et $E= E^{-1,0} \oplus E^{0,-1}$ la d\'ecomposition de Hodge de celui-ci.
\`A l'aide de la trivialisation $\Psi$, on associe \`a toute section $\theta$ de $E$ un champ de vecteurs
sur $A'$ que l'on note \'egalement $\theta$. En particulier, pour tout $\gamma \in \Pi$, on obtient
deux champs de vecteurs $\gamma^{-1,0}$ et $\gamma^{0,-1}$ sur $A'$, o\`u $\gamma = \gamma^{-1,0} +\gamma^{0,-1}$
est la d\'ecomposition de Hodge de $\gamma$ (cf. Partie \ref{champ} pour une expression en coordonn\'ees
de $\gamma^{-1,0}$ et $\gamma^{0,-1}$). Toujours gr\^ace \`a la trivialisation $\Psi$, on construit une d\'ecomposition
du fibr\'e tangent $\cinf$-complexe de $A'$: $T_{\C}A' = (q')^* T_{\C} \h^g \oplus (A' \times \Pi_{\C})$ o\`u
$T_{\C} \h^g$ d\'esigne le fibr\'e tangent $\cinf$-complexe de $\h^g$. La projection canonique
$T_{\C}A' \to A' \times \Pi_{\C}$ donne une forme diff\'erentielle $\nu$ sur $A'$ \`a valeurs dans le fibr\'e $A' \times \Pi_{\C}$
(cf. Partie \ref{formpol} pour une expression de $\nu$ en coordonn\'ees).
Pour $a \in \N^{\geq 1}$ et $\gamma \in \Pi \setminus \{0\}$,
la forme diff\'erentielle $g'_{a,\gamma} \in \Gamma( A', \Omega^{2g-1}_{A'} \otimes \Sym^{a-1} E)$
est d\'efinie par
$$ \displaystyle g'_{a,\gamma} = i_{\gamma^{0,-1}} \left(
       \chi_{\gamma} \times \frac{1}{( \rho(\gamma) - L_{\gamma^{-1,0}})^a} \; vol \times (\gamma^{0,-1})^{a-1}
        \right),$$

\begin{tabular}{lcp{14cm}}
o\`u  & $\bullet$ & $i_{\gamma^{0,-1}}$ est l'op\'erateur de contraction par le champ de vecteurs $\gamma^{0,-1}$, \\
      & $\bullet$ & $\chi_{\gamma}$ est une fonction complexe sur $A'$ d\'efinie par
             $$\chi_{\gamma}([(z,\tau)]) = \exp ( < \gamma \; , pr_1 \circ   \Psi^{-1} ([(z,\tau)]) >_{\C}),$$
	      $pr_1$ d\'esignant la projection canonique $(\Pi_{\R} / \Pi) \times \h^g  \to \Pi_{\R} / \Pi$ et 
	       $< \cdot \; , \cdot >_{\C}$
	      
l'accouplement obtenu en prolongeant par lin\'earit\'e $<  \cdot \; , \cdot >$ (cf. Partie \ref{chi} pour une expression de $\chi_{\gamma}$ en coordonn\'ees), \\
& $\bullet$ & $\rho(\gamma) = <  \gamma^{-1,0} \; , \gamma^{0,-1}  >_{\C}$ (cf. Partie \ref{norme} pour une expression de
      $\rho(\gamma)$ en coordonn\'ees), \\
& $\bullet$ & $L_{\gamma^{-1,0}}$ est la d\'eriv\'ee de Lie associ\'ee au champ de vecteurs $\gamma^{-1,0}$, \\
& $\bullet$ & $vol = (-1)^g (g!)^{-1} \omega^g$, o\`u $\omega := \displaystyle \frac{1}{2} < \nu , \nu  >_{\C}$ (cf. Partie
     \ref{formpol} pour une expression de $\omega$ en coordonn\'ees).\\
     \\
\end{tabular}

La forme $ \displaystyle g'_{a,\gamma}$ est bien d\'efinie car $L_{\gamma^{-1,0}}^k vol = 0$ pour $k > 2g$ (cf. \cite[Prop 3.2.2]{l}). Dans la partie \ref{calctech}, on effectue le calcul en coordonn\'ees de
$$ \displaystyle i_{\gamma^{0,-1}} \left(
       \chi_{\gamma} \times \frac{1}{( \rho(\gamma) - L_{\gamma^{-1,0}})^a} \;\omega^g \right).$$

Dans une derni\`ere sous-partie, on applique les r\'esultats de la Partie 5 de \cite{b}, pour obtenir une expression
explicite des classes d'Eisenstein dans cette situation g\'eom\'etrique (cf. Prop \ref{explicit_Eis}).

\subsection{Polarisation et trivialisation $\cinf$-r\'eelle} \label{trivialisation}

Le calcul des courants de Levin recquiert une polarisation.
On pr\'ecise ici celle que l'on consid\`ere dans la suite. \\

Le pullback de la famille de tores complexes
$q \colon \Lambda'(\mathfrak{a},N) \backslash ( \C^g \times  \h^g) \to \Lambda(\mathfrak{a},N) \backslash \h^g$
par la projection canonique $\h^g \to \Lambda(\mathfrak{a},N) \backslash \h^g$ est
$$
q' \colon A' := \Pi \backslash ( \C^g \times \h^g  ) \to \h^g,
$$
o\`u le quotient de $\C^g \times \h^g$ par le groupe $\Pi := \mathfrak{a}^{\vee} \oplus \mathfrak{a}$ est
celui associ\'e \`a l'action
$$
\begin{array}{ccc}
  (  \mathfrak{a}^{\vee} \oplus \mathfrak{a} ) \times ( \C^g \times \h^g ) & \to &  \C^g \times \h^g  \\
     ( ( a' , a ) , ( z , \tau )) & \mapsto &
     ( (\sigma_k(a')  + \sigma_k(a) \tau_k + z_k )_{1 \leq k \leq g} ,   \tau ) .
\end{array}
$$

Soit $\tau \in \h^g$. On note $\varphi_{\tau}$ le monomorphisme
$\varphi_{\tau} \colon \Pi \to \C^g, \; (a',a) \mapsto ((\sigma_k(a')  + \sigma_k(a) \tau_k)_{1 \leq k \leq g}$.
Le morphisme $\varphi_{\tau} \otimes \R \colon \Pi \otimes \R \to \C^g$ d\'eduit de $\varphi_{\tau}$ par lin\'earit\'e
est un isomorphisme de $\R$-vectoriels. \\

La fibre en $\tau \in \h^g$ de $q'$, not\'ee $A'_{\tau}$, est $\C^g / \varphi_{\tau}(\Pi)$ et
l'application
$$
\begin{array}{cccl}
   <  \cdot \; , \cdot   > \colon & \Pi \wedge \Pi & \to & 2 \pi i \Z \\
           &  (a'_1,a_1) \wedge (a'_2,a_2) & \mapsto &  2 \pi i \; \Tr_L ( a'_2 a_1 - a'_1 a_2)
\end{array}
$$
induit, via $\varphi_{\tau}$, une polarisation principale sur $A'_{\tau}$. \\

On a en outre un
isomorphisme $\cinf$-r\'eel de familles de tores r\'eels au dessus de $\h^g$
$$\begin{array}{cccl}
 \Psi: & ( \Pi_{\R} / \Pi ) \times \h^g & \to &  A' = \Pi \backslash ( \C^g \times \h^g  ).\\
                  & ([\xi],  \tau) & \mapsto &
                  [ (\varphi_{\tau} \otimes \R \; (\xi)  , \tau)] \end{array}
$$
Celui-ci joue un r\^ole fondamental dans la construction des courants de Levin (cf. \cite[Part 2.2]{l} pour la d\'efinition de cette trivialisation dans le cas g\'en\'eral). \\

Il est utile pour la suite d'expliciter l'inverse de $\Psi$. On fixe $(\alpha_1, \dots , \alpha_g)$ une $\Z$-base
de $\mathfrak{a}$ et on note $(\alpha'_1, \dots , \alpha'_g)$ la $\Z$-base de $\mathfrak{a}^{\vee}$ duale de
$(\alpha_1, \dots , \alpha_g)$ relativement à la forme bilinéaire $\mathfrak{a}^{\vee} \times \mathfrak{a} \to \Z$,
$(\alpha' ,\alpha) \mapsto \Tr_L(\alpha' \alpha)$.

\begin{notas} $-$ \\

\begin{itemize}
\item[$\bullet$] Pour $z=(z_1 \dots ,z_g) \in \C^g$, on note $D(z)$ la matrice diagonale
            $g \times g$ dont la diagonale est $(z_1, \dots ,z_g)$ et
             $ \overline{z} := (\overline{z_1}, \dots ,\overline{z_g}).$\\
             
\item[$\bullet$] On introduit les deux matrices $g \times g$ 
         $$ M:= ( \sigma_l (\alpha_k))_{1 \leq k,l \leq g} \quad \mbox{ et } \quad 
              M' := ( \sigma_l(\alpha'_k))_{1 \leq k,l \leq g} .$$
\item[$\bullet$] Pour $\tau = (\tau_1, \dots, \tau_g) \in \h^g$, on pose
              $$ T(\tau) := (t_1,..,t_g) :=
            ((\tau_1-\overline{\tau_1})^{-1}, \dots , (\tau_g-\overline{\tau_g})^{-1}).$$
\end{itemize}
\end{notas}

\`A l'aide des deux $\Z$-bases $(\alpha_1, \dots , \alpha_g)$ et $(\alpha'_1, \dots , \alpha'_g)$,
on effectue un calcul \'el\'ementaire pour montrer que  l'image de 
$[(z, \tau)]  \in A' = \Pi \backslash ( \C^g \times \h^g  )$ par le morphisme $\cinf$-r\'eel
$\Psi^{-1}$ est l'élément de $( \Pi_{\R} / \Pi ) \times \h^g $
$$
\left( 
\left[ 
(z-\overline{z}) D(T(\tau)) M^{-1} \left( \begin{array}{c} \alpha_1 \\ \vdots \\ \alpha_g \end{array} \right) 
+
( \overline{z} D(T(\tau)) D(\tau) (M')^{-1} - z D(T(\tau)) D( \overline{\tau}) (M')^{-1}  ) 
\left( \begin{array}{c} \alpha'_1 \\ \vdots \\ \alpha'_g \end{array} \right) 
\right] ,
\tau \right) 
.$$

\subsection{D\'ecomposition de Hodge} \label{dechodge}

On fixe $\tau \in \h^g$. On note
 $ ( u_{\tau,1}, \dots , u_{\tau,g})$
   (resp. 
$ \overline{u}_{\tau,1} , \dots,
                  \overline{u}_{\tau,g}) $
la base de   duale de $dz:=(dz_1,  \dots ,dz_g)
 \mbox{ (resp. } d\overline{z} :=(d\overline{z_1}, \dots ,d\overline{z_g}) \mbox{)} $
des formes holomorphes (resp. antiholomorphes) de $A'_{\tau}$, la fibre de $q'$ au dessus de $\tau$. \\

L'égalité matricielle suivante nous fournit une expression de 
$ u_{\tau,1}, \dots , u_{\tau,g},  \overline{u}_{\tau,1} , \dots, \overline{u}_{\tau,g}$ 
dans $\Pi_{\C} = H_1(A'_{\tau},\C)$.

$$
\left( 
\begin{array}{c} u_{\tau,1} \\ \vdots \\ u_{\tau,g} \\ \overline{u}_{\tau,1} \\ \vdots \\ \overline{u}_{\tau,g} \end{array}
\right) 
=
\left(
\begin{array}{cc}
D(T(\tau)) M^{-1} & -D(T(\tau)) D(\overline{\tau}) (M')^{-1}  \\
-D(T(\tau)) M^{-1} & D(T(\tau))  D(\tau)  (M')^{-1} \end{array}
\right)
\left( 
\begin{array}{c} \alpha_1 \\ \vdots \\ \alpha_g \\ \alpha'_1 \\ \vdots \\ \alpha'_g \end{array}
\right) .$$
$\;$\\

 On note $ E := E^{-1,0} \oplus E^{0,-1}$ la d\'ecomposition de Hodge du fibr\'e vectoriel analytique au
dessus de $\h^g$ attach\'ee \`a la variation de structures de Hodge pures de poids $-1$
$(R^1 q_* \Z)^{\vee}$. La trivialisation $\cinf$-r\'eelle $\Psi$ introduite dans la partie
\ref{trivialisation} induit une identification entre $\overline{(R^1 q_* \Z)^{\vee}}$ (syst\`eme local
constant) et $\Pi$. \\

Soit $\gamma=(a',a) \in \Pi$. Cet \'el\'ement d\'efinit une section de $q^* E$  au dessus de $A'$ 
$$ A' \to A' \times \Pi_{\C}, \quad  [(z,\tau)] \mapsto ( [(z,\tau)] , \gamma) $$ 
que l'on note \'egalement (abusivement) $\gamma$.
On v\'erifie que la d\'ecomposition de Hodge de $\gamma$ est donn\'ee par:

$$
\begin{array}{cccl}
  \gamma^{-1,0} : &  A' & \to & q^* E^{-1,0} \\
                  & [(z,\tau)] & \mapsto &
          \left( [(z,\tau)] ,
       \displaystyle \sum_{k=1}^g (\sigma_k(a')  + \sigma_k(a) \tau_k)  \;  u_{\tau,k} \right)\\
       \\
\gamma^{0,-1} : & A' & \to & q^* E^{0,-1} \\
                  & [(z,\tau)] & \mapsto &
          \left( [(z,\tau)] ,
       \displaystyle \sum_{k=1}^g ( \sigma_k(a')  + \sigma_k(a) \overline{\tau_k}  )  \; \overline{u}_{\tau,k} \right).\\
\end{array}
$$

\subsection{Champs de vecteurs} \label{champ}

La trivialisation $\Psi$ permet d'associer \`a tout \'el\'ement $\theta$
de $\Pi_{\C}$ un champ de vecteurs $\cinf$ complexe sur $A'$ que l'on note (abusivement) aussi $\theta$.
Pour $\gamma=(a',a) \in \Pi$, les champs de vecteurs associ\'es \`a  $\gamma^{-1,0}$ et $\gamma^{0,-1}$
sont:
$$ \gamma^{-1,0} = \displaystyle \sum_{k=1}^g (\sigma_k(a')  + \sigma_k(a) \tau_k)  \;  \partial_{z_k}
   \quad \mbox{ et } \quad
   \gamma^{0,-1} = \displaystyle \sum_{k=1}^g ( \sigma_k(a')  + \sigma_k(a) \overline{\tau_k}) \;
   \partial_{\overline{z_k}}.$$

\subsection{La forme de polarisation} \label{formpol}

On d\'eduit de la trivialisation $\Psi$ une d\'ecomposition du fibr\'e tangent $\cinf$-complexe de
$A'$: $T_{\C}A' = (q')^* T_{\C} \h^g \oplus (A' \times \Pi_{\C})$. La projection de $TA'$ sur
$(A'\times \Pi_{\C})$ induit une forme diff\'erentielle \`a valeurs dans
$(A' \times \Pi_{\C})$ not\'ee $\nu$.

\begin{notas} $-$ \\

\begin{itemize}
\item[$\bullet$] $D(d\tau)$ (resp. $D(d\overline{\tau})$)
         est la matrice diagonale de formes diff\'erentielles dont la diagonale est
         $(d\tau_1, \dots, d\tau_g)$ (resp. $(d\overline{\tau_1}, \dots, d\overline{\tau_g})$).\\

\item[$\bullet$] $u_k$ (resp. $\overline{u}_k$) est la section de $(A' \times \Pi_{\C})$
au dessus de $A'$
         $ u_k \colon [(z,\tau)] \mapsto u_{\tau,k}$ (resp. $\overline{u}_k
	  \colon [(z,\tau)] \mapsto \overline{u}_{\tau,k}$) pour $k \in \{ 1, \dots , g\}$.
\end{itemize}
\end{notas}

Au prix d'un calcul faisant intervenir l'expression explicite de $\Psi^{-1}$ obtenue pr\'ec\'edemment,
 on obtient l'expression en coordonn\'ees de
$\nu$ suivante:
$$ \nu = ( dz  - (z - \overline{z}) D(T(\tau))
                    D(d\tau) )
         \left( \begin{array}{c} u_1 \\ \vdots \\ u_g \end{array} \right)
          +
          ( d\overline{z}  - (z - \overline{z}  )
           D(T(\tau)) D( d\overline{\tau}  ))
          \left( \begin{array}{c} \overline{u}_1 \\ \vdots \\ \overline{u}_g \end{array} \right).$$

Soit $\tau \in \h^g$. La matrice de Gram de l'accouplement $< \cdot \; , \cdot >_{\C} \colon \Pi_{\C} \wedge \Pi_{\C}  \to \C$ 
(obtenu en \'etendant par lin\'earit\'e $<\cdot \; , \cdot>$) dans la base
$ ( u_{\tau,1} , \dots , u_{\tau,g},   \overline{u}_{\tau,1} , \dots , \overline{u}_{\tau,g} )$ est

$$ ( <  u_{\tau,k} \; ,   \overline{u}_{\tau,l} >_{\C} )_{1 \leq k,l \leq g}
   = - 2 \pi i \; T(\tau). $$

La forme de polarisation $\omega := \displaystyle \frac{1}{2} < \nu , \nu >_{\C}$ est donc:
$$ \displaystyle \omega = - 2 \pi i \sum_{k=1}^g \underset{\eta^1_k}{\underbrace{(
              dz_k - t_k (z_k - \overline{z_k} )d\tau_k)}} \wedge  \underset{\eta^2_k}{\underbrace{( t_k (
              d\overline{z_k} - t_k (z_k - \overline{z_k} ) d \overline{\tau_k})}} .$$

\subsection{Normes} \label{norme}

Soient $\tau \in \h^g$ et $\gamma = (a',a) \in \Pi$. De $   ( <  u_{\tau,k} \; ,   \overline{u}_{\tau,l} >_{\C} )_{1 \leq k,l \leq g}
   = - 2 \pi i \; T(\tau)$ et des expressions obtenues pour $\gamma^{-1,0}$ et
$\gamma^{0,1}$ dans la partie \ref{dechodge}, on d\'eduit que l'image de $\tau$ par  la fonction 
$ \rho(\gamma) =$ \\$ <\gamma^{-1,0},\gamma^{0,1}> \colon \h^g \to \R$ est
$$ \rho(\gamma) (\tau) = - 2 \pi i
   \displaystyle \sum_{k=1}^g t_k |\sigma_k(a') +  \sigma_k(a) \tau_k |^2.$$

\subsection{Les fonctions $\chi_{\gamma}$ ($\gamma \in \Pi$)} \label{chi}

Soit $\tau \in \h^g$. L'isomorphisme de $\R$-vectoriels $\varphi_{\tau} \otimes \R \colon \Pi_{\R} \to \C^g$
induit par passage au quotient un isomorphisme $\cinf$-r\'eel de tores r\'eels
$\overline{\varphi_{\tau} \otimes \R } \colon \Pi_{\R} / \Pi \to A'_{\tau}$. \\

Soit $\gamma = (a',a) \in \Pi$. L'image de $[(z,\tau)] \in A'$ par 
la fonction $\cinf$-r\'eelle
$\chi_{\gamma} \colon A' \to \C$ est donn\'ee par:
$$ \chi_{\gamma} ([(z,\tau)]) = \exp \left( < \gamma ,  \overline{\varphi_{\tau} \otimes \R }^{-1} (z)>_{\C} \right).$$

\subsection{D\'eriv\'ees de Lie successives de la forme $vol$ et effet du produit int\'erieur}
\label{calctech}

Soit $\alpha \in \N^{\geq 1}$ et soit $\gamma = (a',a) \in \Pi \setminus \{0\}$.
On explicite tout d'abord 
$ \displaystyle \frac{1}{(\rho(\gamma) - L_{\gamma^{-1,0}})^{\alpha}} vol $, 
o\`u $vol := (-1)^g (g!)^{-1} \omega^g$. 
On rappelle que $L_{\gamma^{-1,0}}$ d\'esigne la d\'eriv\'ee de Lie du champ de vecteurs $\gamma^{-1,0}$ 
et que
$L_{\gamma^{-1,0}}^k \omega^g =0$ si $k > 2g$ (cf. \cite[Prop. 3.2.2]{l}) . 
\`A l'aide de l'expression de $\omega$ de la Partie \ref{formpol}, on obtient:

$$
\begin{array}{lll}
\displaystyle \frac{1}{(\rho(\gamma) - L_{\gamma^{-1,0}})^{\alpha}} vol & = &  \displaystyle
 \frac{(2 \pi i )^g}{(\rho(\gamma) - L_{\gamma^{-1,0}})^{\alpha}} \;   
 \eta^1_1 \wedge \eta^2_1  \wedge \eta^1_2 \wedge \eta^2_2 \wedge \dots \wedge \eta^1_g \wedge \eta^2_g \\
 \\
 & =  & \displaystyle
 \sum_{n=0}^{2g} \; \frac{(2 \pi i)^g \; C^{\alpha - 1}_{n + \alpha - 1}
 }{\rho(\gamma)^{n + \alpha}}  \; (L_{\gamma^{-1,0}})^n \;
 \eta^1_1 \wedge \eta^2_1  \wedge \eta^1_2 \wedge \eta^2_2 \wedge \dots \wedge \eta^1_g \wedge \eta^2_g.
 
 \end{array}
$$

Pour tout $n \in \{ 0, \dots , 2g \}$, on note $\mathcal{L}_n$ l'ensemble
$$ \{ (L_1, \dots ,L_{2g}) \in \{ \Id ,  L_{\gamma^{-1,0}} \}^{2g} \; | \; \sharp (\{k \in \{ 1 , \dots , 2g \} \; | \;
 L_k =  L_{\gamma^{-1,0}} \}) = n \} .$$
Puisque $(L_{\gamma^{-1,0}})^2 \; \eta_k^1 = (L_{\gamma^{-1,0}})^2 \; \eta_k^2 = 0$ pour tout $k \in \{ 1, \dots , g\}$, , on a:
$$ (L_{\gamma^{-1,0}})^n  \; \eta^1_1 \wedge \eta^2_1 \wedge \dots \wedge \eta^1_g \wedge \eta^2_g 
= \sum_{(L_1, \dots ,L_{2d}) \in \mathcal{L}_n} \;
   L_1 \eta^1_1 \wedge L_2 \eta^2_1 \wedge \dots  \wedge L_{2g-1} \eta^1_g \wedge L_{2g} \eta_g^2 .$$
On vérifie alors que pour tout $ k \in \{ 1, \dots , g\}$:
$$ L_{\gamma^{-1,0}} \; \eta^1_k = -t_k ( \sigma_k(a') + \sigma_k(a) \overline{\tau_k}) d \tau_k
  \quad  \mbox{ et } \quad 
L_{\gamma^{-1,0}} \; \eta^2_k = -t_k^2 ( \sigma_k(a')  + \sigma_k(a) d \overline{\tau_k}) .$$

On a ainsi une expression en coordonnées de $ \displaystyle \frac{1}{(\rho(\gamma) - L_{\gamma^{-1,0}})^{\alpha}} vol $. On calcule maintenant
$i_{\gamma^{0,-1}} \;  \displaystyle \frac{1}{(\rho(\gamma) - L_{\gamma^{-1,0}})^{\alpha}} vol$ ( $i_{\gamma^{0,-1}}$ d\'esigne l'op\'erateur de contraction par le champ de vecteurs $\gamma^{0,-1}$). Compte tenu de ce qui précède, il suffit d'expliciter 
$ i_{\gamma^{0,-1}} L_1 \eta_1^1 \wedge L_2 \eta_1^2 \wedge \dots \wedge L_{2g-1} \eta_g^1 \wedge L_{2g} \eta_g^2$ pour 
tout $n \in \{ 0, \dots , 2g \}$ et $(L_1,..,L_{2g}) \in \mathcal{L}_n$.\\

Soient $n \in \{ 0, \dots , 2g \}$  et $(L_1,..,L_{2g}) \in \mathcal{L}_n$.
$$ \begin{array}{ccl}
  i_{\gamma^{0,-1}} \; L_1 \eta_1^1 \wedge L_2 \eta_1^2 \wedge \dots \wedge L_{2g-1} \eta_g^1 \wedge L_{2g} \eta_g^2& = &
               ( i_{\gamma^{0,-1}} L_1 \eta_1^1 ) \;  L_2 \eta_1^2 \wedge \dots \wedge L_{2g-1} \eta_g^1 \wedge L_{2g}
              \eta_g^2\\
     & - &      ( i_{\gamma^{0,-1}} L_2 \eta_1^2 ) \;  L_1 \eta_1^1 \wedge L_3 \eta_2^1  \wedge \dots  \wedge L_{2g-1} \eta_g^1
              \wedge L_{2g} \eta_g^2\\
     & + &        ... \\
     & - & ( i_{\gamma^{0,-1}} L_{2g} \eta_g^2 )  \; L_1 \eta_1^1 \wedge L_2 \eta_1^2 \wedge \dots  \wedge L_{2g-1} \eta_g^1 .
\end{array}
 $$
De plus, pour tout $k \in \{ 1 , \dots , g \}$:
$$ i_{\gamma^{0-1}} \; \eta^2_k = t_k (  \sigma_k(a') +  \sigma_k(a) \overline{\tau_k} ) \quad  \mbox{ et }\quad
    i_{\gamma^{0-1}} \; \eta^1_k =  i_{\gamma^{0-1}} \; L_{\gamma^{-1,0}} \; \eta^1_k =
    i_{\gamma^{0-1}} \; L_{\gamma^{-1,0}} \; \eta^2_k = 0 .$$

\subsection{D\'etermination des classes d'Eisenstein au niveau topologique}

On fixe $(b',b) \in  (N^{-1} \mathfrak{a}^{\vee} \oplus N^{-1} \mathfrak{a}) \setminus \{ 0 \} $ et
$l \in \N^{>2g}$.\\

De la d\'efinition du morphisme $\PP_{\omega}$ de Levin \cite[Thm 3.4.4]{l} et de la Partie 5 de
\cite{b} et d'un calcul \'el\'ementaire, on d\'eduit que $For(\Eis_{x_{b',b}}^l) \in
H^{2g-1}_{\text{Betti}}( S^0(\C) , \overline{ (\Sym^l  \hh )(g)} )$ co\"incide avec la classe de cohomologie
induite par la forme diff\'erentielle
$$ (-1)^l (l+2g)  \; \sum_{\gamma \in \Pi \setminus \{0\}} \; x_{b',b}^* \;  g'_{l+1,\gamma}
\in \Gamma( \h^g , \Omega_{\h^g}^{2g - 1} \otimes \Sym^l \; \Pi_{\C})
$$
o\`u $  g'_{l+1,\gamma} = \chi_{\gamma} \; i_{\gamma^{0,-1}} \;  \displaystyle \frac{1}{(\rho(\gamma) - L_{\gamma^{-1,0}})^{l+1}} vol \; 
            \otimes (\gamma^{0,-1})^l,$ pour $\gamma \in \Pi \setminus \{0\}$.

Soit $\gamma = (a',a) \in \Pi $ et $k \in \{ 1 , \dots , g \}$.
On explicite maintenant $x_{b',b}^* \; g'_{l+1,\gamma}$ \`a l'aide des calculs effectu\'es pr\'ec\'edemment.\\

On a les relations suivantes
$$\begin{array}{lll}
x_{b',b}^* \; \eta_k^1  = x_{b',b}^* \; \eta_k^2  =  0 &  \quad &
x_{b',b}^* \; L_{\gamma^{-1,0}}  \; \eta_k^1  =  -t_k (\sigma_k(a')   + \sigma_k(a) \overline{\tau_k}) d \tau_k \\
x_{b',b}^* \; L_{\gamma^{-1,0}}  \; \eta_k^2  = -t_k^2 ( \sigma_k(a')   + \sigma_k(a)  \tau_k) d \overline{\tau_k} &  &
x_{b',b}^* \; i_{\gamma^{0-1}} \; \eta_k^2  =  t_k ( \sigma_k(a')   +  \sigma_k(a) \overline{\tau_k}). 
\end{array}
$$
et par suite
$$ (-1)^l (l+2g)   \; x_{b',b}^* \; i_{\gamma^{0,-1}} \;  \displaystyle \frac{1}{(\rho(\gamma) - L_{\gamma^{-1,0}})^{l+1}} vol 
    = 
    \displaystyle
    \frac{(-1)^{l+1} \; (2g +l )!  \; (2 \pi i)^g }{ l! \; \;  \rho(\gamma)^{2g + l } } \; \sum_{k=1}^g  \nu_k $$
avec
$$ \nu_k := t_k^2 (\sigma_k(a')  + \sigma_k(a)  \overline{\tau_k})^2
   \left( \prod\limits_{j \not= k} t_j^3 |\sigma_j(a')  + \sigma_j(a) \tau_j|^2 \right)
  d\tau_1 \wedge d \overline{\tau_1}  \wedge \dots \wedge d\tau_k \wedge \widehat{d \overline{\tau_k}} \wedge \dots
  \wedge d\tau_g \wedge d \overline{\tau_g} ,$$
pour $k \in \{ 1 , \dots , g \}$. En effet,  seule reste la contribution de
$ i_{\gamma^{0,-1}} \; L_{\gamma^{-1,0}}^{2g-1} \; vol .$
On a également
$$ x_{b',b}^* \; \chi_{\gamma} = \exp ( < \gamma , (b',b) >) = \exp( 2 \pi i \; \Tr_L ( a'b - ab') ) .$$
De plus, on a explicité $\gamma^{-1,0}$ relativement à la base
 $ (\overline{u}_{\tau,1}, \dots,  \overline{u}_{\tau,g})  $, mais aussi, la base 
  $ (\overline{u}_{\tau,1}, \dots,  \overline{u}_{\tau,g})  $
relativement à une $\Z$-base de $\mathfrak{a}^{\vee} \oplus \mathfrak{a}$ (cf. Partie \ref{dechodge}). On en déduit l'expression de
$\gamma^{-1,0}$ dans $\C^g \oplus \C^g$, via l'isomorphisme
$\iota \colon ( \mathfrak{a}^{\vee} \oplus \mathfrak{a}) \otimes_{\Z} \C \isom \C^g \oplus \C^g $
 induit par les plongements $(\sigma_k)_{1 \leq k \leq g}$:
$$ \gamma^{-1,0} = ( ( - ( t_k \overline{\tau_k} ( \sigma_k(a') + \sigma_k'(a) \tau_k ) )_{1 \leq k \leq g} ,
                      ( t_k ( \sigma_k(a')  +  \sigma_k'(a) \tau_k ) _{1 \leq k \leq g} ).$$

On rassemble les r\'esultats obtenus dans la Proposition suivante.

\begin{prop} $-$\label{explicit_Eis}
                    La $(g,g-1)$-forme diff\'erentielle sur $\h^g$ \`a valeurs dans
                    $\Sym^l \Pi_{\C}$
	$$ 
	(2 \pi i )^g \; \frac{(-1)^{l+1} \; (2g +l )! }{l!}    \displaystyle \sum_{(a',a) \in \mathfrak{a}^{\vee} \oplus \mathfrak{a}  \setminus \{ 0 \}} \; \sum_{k=1}^g \;
	   \frac{\exp(2 \pi i \; \Tr_L(a'b-b'a))}{\rho(a',a)^{2g + l}} \; f_k(a',a,\tau )\;\mu_k \; \otimes
	   \; h(a',a,\tau)^l  $$
o\`u pour tout $(a',a) \in \mathfrak{a}^{\vee} \oplus \mathfrak{a} $, $\tau \in \h^g$, $k \in \{1 , \dots , g\}$,
on note $t_k  :=(\tau_k- \overline{\tau_k})^{-1}$ et
\begin{center}
\begin{tabular}{l}
$\rho(a',a)   :=
   - 2 \pi i \;
   \displaystyle \sum_{j=1}^g t_j |\sigma_j(a') + \sigma_j(a) \tau_j |^2 $,\\
$f_k(a',a,\tau):=   t_k^2 \; ( \sigma_k(a')  + \sigma_k(a) \overline{\tau_k})^2 \;
    \prod\limits_{j \not= k} t_j^3 |\sigma_j(a')  + \sigma_j(a) \tau_j|^2 $,\\
$\mu_k  := d\tau_1 \wedge d \overline{\tau_1}  \wedge \dots
                 \wedge d\tau_k \wedge \widehat{d \overline{\tau_k}} \wedge \dots
  \wedge d\tau_g \wedge d \overline{\tau_g} \in \Gamma( \h^g,  \Omega^{2g-1}_{\h^g, \C})$, \\
$h(a',a,\tau)  :=   ( - ( t_k \overline{\tau_k} ( \sigma_k(a') + \sigma_k(a) \tau_k ) )_{1 \leq k \leq g} ,
                      ( t_k (  \sigma_k(a') + \sigma_k(a) \tau_k ) _{1 \leq k \leq g} ) \in
		      \C^g \oplus \C^g \overset{\sim}{\underset{\iota}{\leftarrow}}
		       \Pi_{\C}$, \\
\end{tabular}
\end{center}
\noindent induit une classe de cohomologie dans
$H^{2g-1}_{\text{Betti}} ( \Lambda(\mathfrak{a},N) \backslash \h^g , (\overline{(\Sym^l \hh) (d)})_{\C} )$
qui co\"incide avec  $\For ( \Eis_{x_{b',b}}^l )$.
\end{prop}

\section{D\'etermination du r\'esidu en l'$\infty$ des classes d'Eisenstein}

On fixe $b' \in N^{-1} \mathfrak{a}^{-1}$, $b \in (N^{-1} \mathfrak{a}) \setminus \mathfrak{a}$, $\lambda \in \N^{\geq 3}$.\\

On note $pr_{res} : \Sym^{\lambda g} (\C^g \oplus \C^g) \to \C$ la projection qui correspond au niveau complexe
à la projection de
$$ \Sym^{\lambda g} V(\overline{\Q}) \to  \Sym^{\lambda g} V (\overline{\Q})/ W$$
introduite précédemment (cf. Remarque \ref{remproj}).

Le morphisme $\overline{\Res^{\lambda g}_{\infty}} \otimes \C$, obtenu \`a partir de $\overline{\Res^{\lambda g}_{\infty}}$ 
par extension des scalaires de $\Q$ \`a $\C$, est donné  par:
$$\begin{array}{cccl}
\overline{\Res^{\lambda g}_{\infty}} \otimes \C :  &  H^{2g-1}_{\text{Betti}}
( \Lambda(\mathfrak{a}, N) \backslash \h^g, \overline{( \Sym^{\lambda g} \hh (g )})_{\C})
& \to & \C \\
  & \left[ \theta \otimes \left[ v_1 \otimes .. \otimes v_{\lambda g} \right] \right] &
      \mapsto & \left(
    \displaystyle \frac{1}{(2 \pi i)^g} \int_{ \Lambda(\mathfrak{a}, N,\infty) \backslash D_r} \theta \right) \times
       pr_{res} ( \left[ v_1 \otimes .. \otimes v_{\lambda g} \right]).
\end{array} $$
où $D_r := \{ (\tau_1, \dots , \tau_g ) \in \h^g \; | \; \prod\limits_{k=1}^g \Im (\tau_k )= r \}$ pour un nombre réel $r \gg 0$.
Ce morphisme respecte les structures rationnelles sous-jacentes. Cette partie est consacr\'ee au calcul
de
$$
\begin{array}{llll}

 \Res^{\lambda g}_{\infty} (\Eis^{\lambda g}_{x_{b',b}} )& = & \overline{\Res^{\lambda g}_{\infty}} 
 ( For(\Eis^{\lambda g}_{x_{b',b}}) )& \text{(d'après la Proposition \ref{prig})}\\
                      				 & = & \overline{\Res^{\lambda g}_{\infty}} \otimes \C \; ( For(\Eis^{\lambda g}_{x_{b',b}}) )  \\
				   & = &  (-1)^{\lambda g} (\lambda+2)g \displaystyle  \; \sum_{\gamma \in \mathfrak{a}^{\vee} \oplus \mathfrak{a}  \setminus \{0\}} 
				       \overline{\Res^{\lambda g}_{\infty}} \otimes \C ( \; x_{b',b}^* \;  g'_{{\lambda g}+1,\gamma})
				 .\\

\end{array}$$

D'apr\`es la Proposition \ref{explicit_Eis}, le nombre rationnel $\Res^{\lambda g}_{\infty} (\Eis^{\lambda g}_{x_{b',b}} )$ est égal \`a:

$$
\frac{(-1)^{\lambda g +1} \; (( \lambda + 2) g )! }{(\lambda g )!}  \;
\sum_{k=1}^g  \; \;
  \sum_{(a',a) \in \mathfrak{a}^{\vee} \oplus \mathfrak{a}  \setminus \{0\}}  \; \exp( 2 \pi i \; \Tr_L(a'b-b'a))  \; I_{a',a,k} $$
où
$$ I_{a',a,k} := \displaystyle \int_{\Lambda(\mathfrak{a}, N,\infty) \backslash D_r}
                 \frac{ \displaystyle \prod\limits_{j=1}^g  t_j^{\lambda} (\sigma_j(a') + \sigma_j(a) \tau_j)^{\lambda}}{ \displaystyle \left( \sum_{j=0}^g 2
                 \pi i \; t_j|\sigma_j(a') +  \sigma_j(a)\tau_j |^2 \right)^{( \lambda + 2) g }} \; \; \nu_k $$
avec $$ \nu_k  = t_k^2 (\sigma_k(a') + \sigma_k(a) \overline{\tau_k})^2
   \left( \prod\limits_{j \not= k} t_j^3 |\sigma_j(a')  + \sigma_j(a) \tau_j|^2 \right)
  d\tau_1 \wedge d \overline{\tau_1}  \wedge \dots \wedge d\tau_k \wedge \widehat{d \overline{\tau_k}} \wedge \dots
  \wedge d\tau_g \wedge d \overline{\tau_g}, $$
pour $k \in \{ 1 ,  \dots , g \}$.
\vspace{0.5cm}

\begin{itemize}
 \item[$\bullet$ \'Etape 1: ]
Calcul de $I_{1,k} :=\displaystyle  \sum_{a \in \mathfrak{a} \setminus \{0\}}   \sum_{a' \in \mathfrak{a}^{\vee}}
             \exp(2 \pi i \; \Tr_L(a'b-b'a))  \;  I_{a',a,k}$,  pour $k \in \{ 1 ,  \dots , g \}$. \\
             
\begin{itemize}
\item[]
\item[i) ] Dans $I_{a',a,k}$, on peut effectuer ($a \not= 0$) le changement de variables
$$ \tau'_j = \tau_j  + \frac{\sigma_j(a')}{\sigma_j(a)}, \mbox{ pour } 1 \leq j \leq g $$
pour observer que $I_{a',a,k} = I_{0,a,k}$.\\

\item[ii) ]Comme $b \notin \mathfrak{a}$, il existe
$a' \in \mathfrak{a}^{\vee}$ tel que: $$\Tr_L(a'b) \notin \Z.$$

\item[iii) ] De i) et ii), on déduit que $I_{1,k} = 0 $. \\
\end{itemize}

\item[$\bullet$ \'Etape 2: ] Simplification  de
   $ I_{2,k} := \displaystyle  \sum_{a' \in  \mathfrak{a}^{\vee} \setminus \{0\}} \exp(2 \pi i \; \Tr_L(a'b)) \; I_{a',0,k}$,  pour $k \in \{ 1 ,  \dots , g \}$. \\
   
On ne calcule que $I_{2,k}$ que pour $k=1$, la m\'ethode \'etant analogue pour les autres valeurs de $k$
et on peut supposer $r=1$.
De nombreuses simplifications apparaissent lorsque $a=0$ dans $I_{a',a,k}$. On note
$x_k := \Re(\tau_k)$ et  $y_k = \Im(\tau_k)$, pour $k \in \{ 1 ,  \dots , g \}$. On a:

$$ I_{2,k} =   \sum_{a' \in \mathfrak{a}^{\vee} \setminus \{0\}  } \exp(2 \pi i \; \Tr_L(a'b)) \;
    \frac{ {(-1)^{g-1} \; \Norm_L(a')}^{\lambda +2 } }{(2 \pi i)^{(l + 2g)}} \; J_{a'} $$
où $$ J_{a'} =   \int_{\Lambda(\mathfrak{a}, N,\infty) \backslash D_1}\; \frac{1}{ \displaystyle \left( \sum_{j=0}^g
    \frac{\sigma_j(a')^2}{y_j}  \right)^{(\lambda + 2)g}  } \; 
    (dx_1 + i d y_1) \wedge dx_2 \wedge \frac{dy_2}{y_2} \wedge \dots \wedge
    dx_g \wedge \frac{dy_g}{y_g} .$$
De la relation $y_1 \dots y_g=1$ vérifiée par $(\tau_1, \dots, \tau_g) \in D_1$, on déduit:
$$ 
\begin{array}{ccl}
J_{a'}  & = &  \displaystyle  \int_{\Lambda(\mathfrak{a}, N,\infty) \backslash D_1} \; \frac{1}{\left( \displaystyle \sum_{j=0}^g
    \frac{\sigma_j(a')^2}{y_j}  \right)^{(\lambda + 2)g}  } \; 
    dx_1 \wedge dx_2 \wedge \frac{dy_2}{y_2} \wedge \dots \wedge
    dx_g \wedge \frac{dy_g}{y_g} \\ 
    \\
%
%
%
%
 &  = &  \displaystyle  \sum_{a' \in ( \mathfrak{a}^{\vee} \setminus \{0\} ) / U_{L,N}}
        \frac{ (-1)^{g-1}  \; \exp(2 \pi i \; \Tr_L(a'b)) \; \text{N}_L(a')^{\lambda +2 } }{(2 \pi i)^{(\lambda + 2)g}}  \\
  & \times &  \displaystyle  \sum_{\varepsilon \in U_{L,N}} 
 \int_{  \Lambda(\mathfrak{a}, N,\infty) \backslash D_1       } \; \frac{1}{\left( \displaystyle  \sum_{j=0}^g
    \frac{ \sigma_j(\varepsilon a')^2}{y_j}  \right)^{(\lambda + 2)g}  } \; dx_1 \wedge dx_2 \wedge \frac{dy_2}{y_2} \wedge \dots \wedge
    dx_g \wedge \frac{dy_g}{y_g} \\
    \\
 &  = & \displaystyle  \sum_{a' \in ( \mathfrak{a}^{\vee} \setminus \{0\} ) / U_{L,N} }
        \frac{ (-1)^{g-1}  \; \exp(2 \pi i \; \Tr_L(cb')) \; \text{N}_L(a')^{\lambda +2 } }{(2 \pi i)^{(\lambda + 2)g}} \; 
        vol(N \mathfrak{d}_L^{-1} \mathfrak{a}^{-2}) \; K_{a'} \\
 \end{array}
$$
où $$ K_{a'} =
\displaystyle  
      \int_{(\R^{>0})^{g-1}} \; \frac{1}{\left( \displaystyle  \sum_{j=0}^g
    \frac{\sigma_j(a')^2}{y_j}  \right)^{(\lambda + 2)g}  }
         \frac{dy_2}{y_2} \dots  \frac{dy_g}{y_g},
 $$    
pour  $a' \in ( \mathfrak{a}^{\vee} \setminus \{0\} ) / U_{L,N}.$
Pour établir la dernière égalité, on utilise le Lemme $2.10_1$ de [F] et l'injectivité de l'homomorphisme de groupes
$U_{L,N} \to U_{L,N}$, $\varepsilon \mapsto \varepsilon^2$.\\

\item[$\bullet$ \'Etape 3: ] Calcul de $K_{a'}$, pour  $a' \in ( \mathfrak{a}^{\vee} \setminus \{0\} ) / U_{L,N}.$\\

On calcule en fait $ \Gamma((\lambda + 2) g) K_{a'} $:
$$ \Gamma((\lambda + 2) g) K_{a'}  =
\int_{(\R^{>0})^g} \; u^{(\lambda + 2) g} \exp 
\left( -u \left(  y_2 \dots y_g \sigma_1(a')^2 + \frac{\sigma_2(a')^2}{y_2} + \dots +\frac{\sigma_g(a')^2}{y_g} \right) \right)
     \frac{du}{u} \frac{dy_2}{y_2} \dots \frac{dy_g}{y_g} .$$
On effectue ensuite le changement de variables:
$$ \begin{array}{l}
u_1 = u y_2 \dots y_g \\
u_2 = u / y_2 \\
\quad \quad \vdots \\
u_g = u / y_g
\end{array}$$
et on trouve:
$$ \begin{array}{ccl}
\Gamma((\lambda + 2) g) K_{a'}  & = &
\displaystyle  \left( \int_{\R^{>0}} \; u_1^{\lambda + 2} \exp(-\sigma_1(a')^2u_1) \frac{du_1}{u_1} \right)  \dots
                         
\left( \int_{\R^{>0}} \; u_g^{\lambda + 2} \exp(-\sigma_g(a')^2u_g) \frac{du_g}{u_g} \right) \\
\\
  & = & \displaystyle \frac{\Gamma(\lambda + 2)^g}{\text{N}_L(a')^{2(\lambda + 2)}}.
 \end{array}  
$$

Ainsi $K_{a'} = \displaystyle \frac{ ((\lambda+1)!)^g}{((\lambda + 2) g - 1)! \; \text{N}_L(a')^{2(\lambda+2)}}$.
\end{itemize}

\begin{nota} $-$ \label{fl}
Pour $\beta \in (N^{-1} \mathfrak{a}) \setminus  \mathfrak{a}$, $s \in \C$ tel que $\Re(s) > 1$,
on d\'efinit $\mathfrak{L}(\mathfrak{a},N,\beta,s)$ par
$$ \mathfrak{L}(\mathfrak{a} ,N,\beta,s) := \displaystyle
                \sum_{ a' \in (\mathfrak{a}^{\vee} \setminus \{0\}
		)/ U_{L,N}}
                \frac{\exp(2 \pi i \; \Tr_L(a' \beta ))}{\text{N}_L(a')^{s}}.$$
\end{nota}

Du calcul que l'on vient d'achever,  on d\'eduit le Théorème suivant.

\begin{thm} $-$ Soient $\lambda \in \N^{\geq 3}$, \label{R1}
                 $b' \in N^{-1} \mathfrak{a}^{\vee}$ et
		   $b \in (N^{-1} \mathfrak{a}) \setminus  \mathfrak{a} $. 

$$ \Res^{\lambda g}_{\infty}(\Eis^{\lambda g}_{x_{b',b}}) =
\displaystyle  \frac{ (-1)^{(\lambda + 1) g} \; ( ( \lambda + 1) ! )^g  \; (\lambda + 2) g^2 \;  N^g }{(\lambda g)! \; \mathcal{N}_L(\mathfrak{a})^2 } \; 
\frac{1}{\sqrt{d_L} \; (2 \pi i)^{(\lambda + 2)g} } \;    \mathfrak{L}(\mathfrak{a},N,b, \lambda + 2)$$
\vspace{0.5cm}
\end{thm}

Ayant pris soin de respecter les structures rationnelles tout au long du calcul,
on  déduit de ce r\'esultat le Théorème de Klingen-Siegel.

\begin{coro}{}$-$ \label{R2}
Pour tout $ \lambda \in \N^{\geq 6}$ pair et $b \in (N^{-1} \mathfrak{a}) \setminus  \mathfrak{a}$:
         $$ \mathfrak{L}(\mathfrak{a},N,b, \lambda) \;  (2 \pi i)^{-\lambda g} \; \sqrt{d_L} \in \Q.$$
\end{coro}

\begin{coro} $-$  Soient $ \lambda \in \N^{\geq 6}$ pair, \label{R3}
$b' \in N^{-1} \mathfrak{a}^{\vee}$
et $b \in (N^{-1}  \mathfrak{a}) \setminus  \mathfrak{a}   $ tel que les idéaux entiers
$N \OO_L $ et $ N b \mathfrak{a}^{-1}$ sont copremiers. Si $g \geq 2$, alors on a:
$$  \Eis^{\lambda g}_{x_{b',b}} \not= 0.$$
\end{coro}

\begin{dem} On prouve que   $\Res^{\lambda g}_{\infty}(\Eis^{\lambda g}_{x_{b',b}}) \not= 0$ \`a l'aide du Th\'eor\`eme
\ref{R1}.
On pose  $\mathfrak{f}:= N \OO_L$ et $\mathfrak{b}:=  N b \mathfrak{a}^{-1}$. On introduit:\\

\begin{itemize}

\item[$\bullet$] l'ensemble $\mathcal{E}( \mathfrak{b}, \mathfrak{f} )$ des idéaux entiers $\mathfrak{g}$
          premiers à $\mathfrak{f}$ pour lesquels il existe $\mu \in L$ totalement positif et congru
	   	  à 1 modulo $\mathfrak{f}\mathfrak{b}^{-1}$ tel que:
 	  $$ \mathfrak{g} \mathfrak{b}^{-1} = (\mu),$$
\item[$\bullet$] la fonction holomorphe 
           $\zeta(\mathfrak{b}, \mathfrak{f}, \cdot )$ définie, pour $s \in \C$ tel que $\Re(s)>1$, par:
 	  $$ \zeta(\mathfrak{b}, \mathfrak{f}, s) := \sum_{\mathfrak{g} \in  \mathcal{E}(\mathfrak{b},
 	  \mathfrak{f} )} \quad  N_L(\mathfrak{g})^{-s},$$
	  
\item[$\bullet$] 	$U^+_{L,N}$ le sous-groupe de $U_{L,N}$ formé par les éléments de $U_{L,N}$ totalement positifs, \\
	  
\item[$\bullet$] la fonction holomorphe $ \mathfrak{L}^+(\mathfrak{a} ,N,b, \cdot)$ définie, pour $s \in \C$ tel que $\Re(s)>1$, par:
$$ \mathfrak{L}^+(\mathfrak{a} ,N,b,s) := \displaystyle
                \sum_{ a' \in (\mathfrak{a}^{\vee} \setminus \{0\}
		)/ U^+_{L,N}}
                \frac{\exp(2 \pi i \; \Tr_L (   a' b))}{|\text{N}_L(a')|^{s}}.$$ 	  
 \end{itemize}

 La fonction $\zeta(\mathfrak{b}, \mathfrak{f}, \cdot)$ a un prolongement holomorphe sur $\C -\{1\}$.
 D'après une équation fonctionnelle pour 
 $\mathfrak{L}^+(\mathcal{D}_L^{-1},N,b,\cdot)$ établie par Siegel (cf. \cite[Formule (10)]{si}), on a:

 $$   [ U_{L,N} :  U^+_{L,N}   ] \; \mathfrak{L}(\mathcal{D}_L^{-1},N,b',\lambda + 2)  = 
        \mathfrak{L}^+(\mathcal{D}_L^{-1},N,b',\lambda + 2)  \underset{\R^{\times}}{\sim}
       \zeta(\mathfrak{b}, \mathfrak{f},  -\lambda - 1),      $$
 et, comme  $\lambda$ est pair, $\zeta(\mathfrak{b}, \mathfrak{f},  -\lambda - 1)$ est non nul (cf. \cite[Theorem VII-5.9]{ne}).

\end{dem}

\end{document}